\documentclass[11pt]{amsart}
\usepackage{amsbsy}
\usepackage{amsmath}
\usepackage{amsbsy}
\usepackage{graphicx}
\usepackage{enumerate}
\usepackage[cp1250]{inputenc}
\usepackage{indentfirst}
\usepackage{latexsym}
\usepackage{amssymb}
\usepackage{amsfonts}
\usepackage{color}
\usepackage{dsfont}
\usepackage{mathtools}
\usepackage{amsthm}

\usepackage{tikz}
\usetikzlibrary{arrows}
\usepackage{xcolor}
\usepackage{xypic}
\usepackage{cite}
\usepackage{float}
\usepackage{filecontents}
\usepackage{graphicx}
\usepackage{epstopdf}
\usepackage{epsfig}

\definecolor{darkgreen}{rgb}{0,.7,.3}

\theoremstyle{definition}
\newtheorem{definition}{Definition}[section]
\newtheorem{remark}[definition]{Remark}

\theoremstyle{plain}
\newtheorem{lemma}[definition]{Lemma}
\newtheorem{theorem}[definition] {Theorem}
\newtheorem{proposition}[definition] {Proposition}

\newtheorem{corollary}[definition]{Corollary}

\title{
Exponential equations in acylindrically hyperbolic groups}
\author{Agnieszka Bier}
\address{Department of Applied Mathematics, Silesian Univesity of Technology,
ul. Kaszubska 23, 44 - 101 Gliwice, Poland}
\email{agnieszka.bier@polsl.pl}

\author{Oleg Bogopolski}
\address{{Sobolev Institute of Mathematics of Siberian Branch of Russian Academy
of Sciences, Novosibirsk, Russia}\newline
{and D\"{u}sseldorf University, Germany}}
\email{Oleg$\_$Bogopolski@yahoo.com}

\begin{document}


\maketitle

\begin{abstract} Let $G$ be an acylindrically hyperbolic group and $E$ an exponential equation over $G$.
We show that if $E$ is solvable in $G$, then there exists a solution whose components, corresponding to loxodromic elements, can be linearly estimated in terms of lengths of the coefficients of $E$.
We give a more precise answer in the case where $G$ is a relatively hyperbolic group. Under some assumption of general character, the solvability and the search problems for exponential equations over $G$ can be reduced to its peripheral subgroups.
\end{abstract}

\section{Introduction}

In 2015, Myasnikov, Nikolaev and Ushakov  initiated the study of exponential equations in groups~\cite{MNU}
which has become a topic of intensive investigations on the edge of group theory and complexity theory~\cite{Bogo_Ivanov, KLZ, Dudkin, LZ_1, LZ_2, Frenkel, GKLZ, Figelius, MT,
Lohrey_1,Lohrey_2}.
The results obtained in~\cite{MNU} for hyperbolic groups motivated us to investigate the decidability
of exponential equations in the wider classes of relatively hyperbolic and acylindrically hyperbolic groups.

\begin{definition}
An {\it exponential equation} over a group $G$ is an equation of the form
$$
a_1g_1^{x_1}a_2g_2^{x_2}\dots a_ng_n^{x_n}=1,\eqno{(1.1)}
$$
where $a_1,g_1,\dots, a_n,g_n$ are elements from $G$ and $x_1,\dots,x_n$ are variables (which take values in $\mathbb{Z}$).
A tuple $(k_1,\dots,k_n)$ of integers is called a {\it solution} of this equation if $a_1g_1^{k_1}a_2g_2^{k_2}\dots a_ng_n^{k_n}=1$ in $G$.
\end{definition}

The first main theorem of this paper, Theorem A, is formulated and proved in Section~7. Here we give a simplified
version of this theorem, Theorem A$'$.
It says that if $G$ is an acylindrically hyperbolic group and the above equation is solvable,
then there exists a solution $(k_1,\dots,k_n)$ such that $|k_j|$ corresponding to loxodromic $g_j$
can be linearly bounded in terms of the lengths of the coefficients of this equation.

\medskip


\noindent
{\bf Theorem A$'$.} {\rm (see Theorem~A)}
{\it Let $G$ be an acylindrically hyperbolic group with respect to a generating set $X$.
Then there exists a constant $M>1$ such that for any exponential equation
$$
a_1g_1^{x_1}a_2g_2^{x_2}\dots a_ng_n^{x_n}=1
$$
with constants $a_1,g_1,\dots,a_n,g_n$ from $G$ and variables $x_1,\dots,x_n$, if this equation is solvable over $\mathbb{Z}$, then there exists a solution $(k_1,\dots,k_n)$ with
$$
|k_j|\leqslant \Bigl(n^2+\overset{n}{\underset{i=1}\sum}\,|a_i|_X+\overset{n}{\underset{i=1}{\sum}}\, |g_i|_X\Bigr)\cdot M
$$
for all $j$ corresponding to loxodromic $g_j$.

If, additionally, $G$ is generated by a finite subset $Y$,
then the above estimation remains valid if we replace there $X$ by $Y$ and $M$ by $M \underset{y\in Y}\sup |y|_X$.}

\begin{remark}
The main result of the paper~\cite{MNU} of Myasnikov, Nikolaev and Ushakov says that if $G$ is a hyperbolic group
with a finite generating set $X$,
then there exists a polynomial $p_n(x)$ such that for any exponential equation of the form (1.1), if this equation is solvable then there exists a solution $(k_1,\dots,k_n)$ with
$$|k_j|\leqslant p_n\Bigl(\overset{n}{\underset{i=1}{\sum}}\,|a_i|_X+\overset{n}{\underset{i=1}{\sum}} |g_i|_X\Bigr)$$ for $j=1,\dots ,n$.

We consider the more general case where $G$ is an acylindrically hyperbolic group.
This case is more difficult since $G$ is not necessarily finitely generated in general. Moreover, even if $G$ is  finitely generated, it can happen that $G$ does not act acylindrically on a locally finite graph.

Theorem~A$'$ restricted to the case where $G$ is hyperbolic and $X$ is finite implies that the above polynomial $p_n(x)$ can be taken to be linear. Indeed, all non-loxodromic elements of $G$ have finite orders in this case, and these orders can be bounded from above by a universal constant depending only on $|X|$ and $\delta$, where $\delta$ is the hyperbolicity constant of $G$ with respect to $X$ (see~\cite{Bog_2} or~\cite{BH}).
Theorem A in Section~7 gives further improvements.
\end{remark}

\begin{remark}
The following example shows that the word {\it loxodromic} in the formulation of Theorem A$'$ cannot be omitted even in the case of finitely presented relatively hyperbolic groups.

\medskip

{\it Example.} Let $H$ be a finitely presented group containing the group of rational numbers~$\mathbb{Q}$.
Such group can be constructed using Higman's embedding theorem.
Then the free product $G=H\ast F_2$ is finitely presented and relatively hyperbolic with
respect to the subgroups $H$ and $F_2$, and the elements of $H$ and $F_2$
and their conjugates are elliptic with respect to the generating set $X=H\cup F_2$.
We consider the rational numbers $a=-1$ and $b_i=\frac{1}{i}$ for $i\geqslant 1$ as elements of $G$.
For each $i\in \mathbb{N}$ the exponential equation $ab_i^{x}=1$ has a unique solution (namely $i$),
and the sum of lengths of its coefficients is $|a|_X+|b_i|_X=2$. Thus, there does not exist a function $f$
such that, for all $i$, the solution of $ab_i^{x}=1$ is bounded from above by $f(|a|_X|+|b_i|_X)$.
\end{remark}

\medskip

Theorem B (see Subsection 8.2) deals with certain exponential equations in groups with hyperbolically embedded subgroups; we use it to deduce Theorem~C.
Theorem~C, comparing with Theorem A, gives more information in the case where $G$ is a finitely generated relatively hyperbolic group. It says that for any exponential equation $E$ over $G$, there exists a finite disjunction $\Phi$ of finite systems of exponential equations over peripheral subgroups of $G$
such that $E$ is solvable if and only if $\Phi$ is solvable. If some additional data are known,
one can find such $\Phi$ algorithmically. Moreover, having a solution of $\Phi$, one can find a solution of $E$.


\medskip

\noindent
{\bf Theorem C.}
{\it Let $G$ be a group relatively hyperbolic with respect to a finite collection of
subgroups $\{H_1,\dots,H_m\}$.
Suppose that $G$ is finitely generated, each subgroup $H_i$ is given by a recursive presentation and has solvable word problem,
$G$ is given by a finite relative presentation $\mathcal{P}=\langle X\,|\, \mathcal{R}\rangle$  with respect to  $\{H_1,\dots,H_m\}$, where $X$ is a finite set generating $G$, and that the hyperbolicity constant $\delta$ of the Cayley graph $\Gamma(G,X\cup \mathcal{H})$ is known, $\mathcal{H}=\overset{m}{\underset{i=1}{\bigsqcup}}H_i$.

Then there exists an algorithm which for any exponential equation $E$ over $G$ finds a finite disjunction $\Phi$ of finite systems of equations,
$$
\Phi:=\overset{k}{\underset{i=1}{\bigvee}}\overset{\ell_i}{\underset{j=1}{\bigwedge}} E_{ij},
$$
such that
\begin{enumerate}
\item[{\rm (1)}] each $E_{ij}$ is  an exponential equation over $H_{\lambda}$ for some $\lambda\in \{1,\dots,m\}$ or a trivial equation of kind $g_{ij}=1$, where $g_{ij}$ is an element of $G$,

\item[{\rm (2)}] for any $i=1,\dots, k$, the sets of variables of
$E_{i,j_1}$ and $E_{i,j_2}$ are disjoint if $j_1\neq j_2$,

\item[{\rm (3)}] $E$ is solvable if and only if $\Phi$ is solvable.\\
Moreover, any solution of $\Phi$ can be algorithmically extended to a solution of~$E$.

\end{enumerate}
}

\medskip

In the proof of Theorem A, which is a stronger version of Theorem A$'$, we use the following theorem
about conjugator lengths in acylindrically hyperbolic groups. This theorem seems to be interesting for its own sake.




\begin{theorem}\label{third_main theorem} Let $G$ be an acylindrically hyperbolic group
with respect to a generating set $X$. Let $\delta$ be the hyperbolicity constant of the Cayley graph $\Gamma(G,X)$
and let $N$ be the function from Definition~\ref{acyl_action}. Then there exists a universal constant $C$
such that for any two conjugate elements $h_1,h_2\in G$ of (possibly infinite) order larger than $N(8\delta+1)$,
there exists $g\in G$ such that $h_2=gh_1g^{-1}$ and $|g|_X\leqslant C(|h_1|_X+|h_2|_X)$.
\end{theorem}


\begin{remark} In~\cite{MNU}, the problem about decidability of equations (1.1) in integer numbers is called the Integer Knapsack Problem (IKP) for the group $G$. If we are looking for nonnegative integer solutions, the problem is called the Knapsack Problem (KP). Clearly, the decidability of (IKP) for $G$ implies the decidability of (KP) for $G$.
To our best knowledge the answer to the following problem is unknown.

\medskip

\noindent
{\bf Problem.} Does there exist a finitely presented group $G$ for which the Integer Knapsack Problem
is decidable and the Knapsack Problem is undecidable?

\medskip

\end{remark}

\section{Definitions and preliminary statements}

We introduce general notation and recall some relevant definitions and statements
from the papers~\cite{Bog_1,DOG,Osin_1}.
In this paper, all actions of groups on metric spaces are assumed to be isometric.

\subsection{General notation} All generating sets considered in this paper are assumed to be symmetric, i.e., closed under taking inverse elements.
Let $G$ be a group generated by a subset $X$. For $g\in G$ let $|g|_X$ be the length of a shortest word in $X$ representing $g$. The corresponding metric
on $G$ is denoted by ${d}_X$ (or by ${d}$ if $X$ is clear from the context); thus ${d}_X(a,b)=|a^{-1}b|_X$. The right Cayley graph of $G$ with respect to $X$ is denoted by $\Gamma(G,X)$.
By a path $p$ in the Cayley graph we mean a combinatorial path; the initial and the terminal vertices of $p$ are denoted by $p_{-}$ and $p_{+}$, respectively.
The length of $p$ is denoted by $\ell(p)$. The label of $p$ (which is a word in the alphabet $X$)
is denoted by ${\bold{ Lab}}(p)$.


Recall that a path $p$ in $\Gamma(G,X)$ is called ($\varkappa,\varepsilon)$-{\it quasi-geodesic}, where $\varkappa\geqslant 1$,
$\varepsilon\geqslant 0$, if ${d}(q_{-},q_{+})\geqslant \frac{1}{\varkappa}\ell(q)-\varepsilon$ for any subpath $q$ of $p$.

\subsection{Hyperbolic spaces}





A geodesic metric space $\frak{X}$ is called {\it $\delta$-hyperbolic} if each side of any geodesic triangle $\Delta$ in $\frak{X}$ lies in the $\delta$-neighborhood of the union of the other two sides of $\Delta$.
We will use the following standard facts about hyperbolic spaces.

\begin{lemma}\label{polygon}
Let $\frak{X}$ be a $\delta$-hyperbolic space. Suppose that $R$ is a geodesic $n$-gon in $\frak{X}$.
Then any side of $R$ is at distance at most $(n-2)\delta$
from the union of the other sides of $R$.
\end{lemma}

\begin{lemma}\label{Hausdorff} {\rm (see~\cite[Chapter III.H, Theorem 1.7]{BH})} For all $\delta\geqslant 0$, $\varkappa\geqslant 1$, $\epsilon\geqslant 0$, there exists a constant
$\mu=\mu(\delta,\varkappa,\epsilon)> 0$ with the following property:

If $\frak{X}$ is a $\delta$-hyperbolic space, $p$ is a $(\varkappa,\epsilon)$-quasi-geodesic in $\frak{X}$,
and $[x,y]$ is a geodesic segment joining the endpoints of $p$, then the Hausdorff distance between $[x,y]$
and the image of $p$ is at most $\mu$.
\end{lemma}

The following corollary is a slight generalization of the previous one.

\begin{corollary}\label{close_qg} 
Let $\frak{X}$ be a $\delta$-hyperbolic space, let $p$ and $q$ be $(\varkappa,\epsilon)$-quasi-geodesics in $\frak{X}$ with $\max\{d(p_{-},q_{-}),d(p_{+},q_{+})\}\leqslant r$. Then the
Hausdorff distance between the images of $p$ and $q$ is at most $r+2\delta+2\mu$,
where $\mu=\mu(\delta,\varkappa,\epsilon)$ is the constant from Lemma~\ref{Hausdorff}.
\end{corollary}

\begin{lemma}\label{HausdorffQuasiGeod}
{\rm (see~\cite[Chapitre 3, Th$\acute{\text{\rm e}}$or$\grave{\text{\rm e}}$me 3.1]{CDP})}
For all $\delta\geqslant 0$, $\varkappa\geqslant 1$, $\varepsilon\geqslant 0$, there exists a constant
$\mu=\mu(\delta,\varkappa,\varepsilon)> 0$ with the following property:

If  $\frak{X}$ is a $\delta$-hyperbolic space and $p$ and $q$ are infinite $(\varkappa,\varepsilon)$-quasi-geodesics in  $\frak{X}$ with the same limit points on the Gromov boundary $\partial \frak{X}$, then the Hausdorff distance between $p$ and $q$
is at most $\mu(\delta,\varkappa,\varepsilon)$.
\end{lemma}


The following lemma enables to estimate a displacement of a point on a segment of a hyperbolic space (under the action of an isometry)
via displacements of the endpoints of this segment.

\begin{lemma}\label{equidist}
Let $G$ be a group acting on a $\delta$-hyperbolic space $\frak{X}$. Let $g\in G$ be an element and $[A,B]$ a geodesic in $\frak{X}$.
Suppose that $C$ is a point on $[A,B]$ such that
$d(A,C)>d(A,gA)+2\delta$ and $d(C,B)>d(B,gB)+2\delta$.
Then
$$
d(C,gC)\leqslant 4\delta+\min\{d(A,gA),d(B,gB)\}.
$$
\end{lemma}

\medskip

{\it Proof.}
By assumptions  the distance from $C$ to $[A,gA]\cup [B,gB]$ is larger than $2\delta$.
By Lemma~\ref{polygon}, there exists
a point $D\in [gA,gB]$ such that $d(C,D)\leqslant 2\delta$.
Then
$$
\begin{array}{ll}
d(C,gC) & \leqslant d(C,D)+d(D,gC)\vspace*{2mm}\\
& = d(C,D)+|d(D,gA)-d(gC,gA)|\vspace*{2mm}\\
& = d(C,D)+|d(D,gA)-d(C,A)|\vspace*{2mm}\\
& \leqslant d(C,D)+d(C,D)+d(A,gA)\vspace*{2mm}\\
& \leqslant 4\delta+d(A,gA).
\end{array}
$$
Analogously, we obtain $d(C,gC)\leqslant 4\delta+d(B,gB)$. \hfill $\Box$

\medskip

Without loss of generality, we may assume that $\delta$ is integer.

\medskip

The following lemma will be used in the proof of the elliptic case of Theorem~\ref{third_main theorem}.

\begin{lemma}\label{quasi-geod_2} {\rm (see~\cite[Lemma 4.8]{Bog_1})}
For every $\delta\geqslant 0$, there exists $\varepsilon_1=\varepsilon_1(\delta)\geqslant 0$ such that the following holds.
Suppose that the Cayley graph of a group $G$ with respect to a generating set $X$ is $\delta$-hyperbolic
for some integer $\delta\geqslant 0$.
Let $a,b\in G$ be conjugate elements satisfying $|a|_X\geqslant |b|_X+4\delta+2$.
Then there exist $x,y\in G$ with the following properties:

\begin{enumerate}
\item[{\rm (1)}] $a=x^{-1}yx$;
\vspace*{1mm}
\item[{\rm (2)}] $|y|_X\in \{|b|_X+4\delta+1,|b|_X+4\delta+2\}$;
\vspace*{1mm}
\item[{\rm (3)}] any path $q_0q_1q_2$ in $\Gamma(G,X)$, where $q_0,q_1,q_2$ are geodesics with labels representing $x^{-1},y,x$, is
a $(1,\varepsilon_1)$-quasi-geodesic.
\end{enumerate}
\end{lemma}

\subsection{Two equivalent definitions of acylindrically hyperbolic groups}

\begin{definition}\label{acyl_action} {\rm (see~\cite{Bowditch} and Introduction in~\cite{Osin_1})
An action of a group $G$ on a metric space $S$ is called
{\it acylindrical}
if for every $\varepsilon>0$ there exist $R,N>0$ such that for every two points $x,y$ with $d(x,y)\geqslant R$,
there are at most $N$ elements $g\in G$ satisfying
$$
d(x,gx)\leqslant \varepsilon\hspace*{2mm}{\text{\rm and}}\hspace*{2mm} d(y,gy)\leqslant \varepsilon.
$$
}
\end{definition}

Given a generating set $X$ of a group $G$, we say that the Cayley graph $\Gamma(G,X)$ is
{\it acylindrical} if the left action of $G$ on $\Gamma(G,X)$ is acylindrical.
For Cayley graphs, the acylindricity condition can be rewritten as follows:
for every $\varepsilon>0$ there exist $R,N>0$ such that for any $g\in G$ of length $|g|_X\geqslant R$
we have
$$
\bigl|\{f\in G\,|\, |f|_X\leqslant \varepsilon,\hspace*{2mm} |g^{-1}fg|_X\leqslant \varepsilon \}\bigr|\leqslant N.
$$

Recall that an action of a group $G$ on a hyperbolic space $S$ is called {\it elementary} if the limit set
of $G$ on the Gromov boundary $\partial S$ contains at most 2 points.

\begin{definition}\label{Definition_of_Osin} {\rm (see~\cite[Definition 1.3]{Osin_1})
A group $G$ is called {\it acylindrically hyperbolic} if it satisfies one of the following equivalent
conditions:

\begin{enumerate}
\item[(${\rm AH}_1$)] There exists a generating set $X$ of $G$ such that the corresponding Cayley graph $\Gamma(G,X)$
is hyperbolic, $|\partial \Gamma (G,X)|>2$, and the natural action of $G$ on $\Gamma(G,X)$ is acylindrical.

\medskip

\item[(${\rm AH}_2$)] $G$ admits a non-elementary acylindrical action on a hyperbolic space.
\end{enumerate}
}
\end{definition}

In the case (AH$_1$), we also write that $G$ is {\it acylindrically hyperbolic with respect to $X$}.

\medskip

\subsection{Elliptic and loxodromic elements in acylindrically hyperbolic groups}

The following definition is standard.

\begin{definition}
{\rm
Given a group $G$ acting on a metric space $S$, an element $g\in G$ is called {\it elliptic}
if some (equivalently, any) orbit of $g$ is bounded, and {\it loxodromic} if the map
$\mathbb{Z}\rightarrow S$ defined by
$n\mapsto g^nx$ is a quasi-isometric embedding for some (equivalently, any) $x\in S$. That is,
for $x\in S$, there exist $\varkappa\geqslant 1$ and $\varepsilon\geqslant 0$ such that for any $n,m\in \mathbb{Z}$ we have
$$
d(g^nx,g^mx)\geqslant \frac{1}{\varkappa} |n-m|-\varepsilon.
$$

Let $X$ be a generating set of $G$.
We say that $g\in G$ is {\it elliptic (respectively loxodromic) with respect to $X$} if $g$ is elliptic (respectively loxodromic) for the canonical left action of $G$ on the Cayley graph $\Gamma(G,X)$.
If $X$ is clear from a context, we omit the words ``with respect to $X$''.

The set of all elliptic
(respectively loxodromic) elements of $G$ with respect to $X$ is denoted by ${\rm Ell}(G,X)$ (respectively by ${\rm Lox}(G,X))$.
}
\end{definition}

Note that for groups acting on geodesic hyperbolic spaces,
there is only one additional isometry type of an element- parabolic
(see e.g.~\cite[Chapitre 9, Th$\acute{\text{e}}$or$\grave{\text{e}}$me~2.1]{CDP}).



Bowditch~\cite[Lemma 2.2]{Bowditch} proved that every element of a group acting acylindrically on a hyperbolic space is either elliptic or loxodromic (see a more general statement in~\cite[Theorem 1.1]{Osin_1}).

Recall that any loxodromic element $g$ in an acylindrically hyperbolic group $G$ is contained in a
unique maximal virtually cyclic subgroup~\cite[Lemma 6.5]{DOG}. This subgroup, denoted by $E_G(g)$, is called the {\it elementary subgroup associated with $g$}; it can be described as follows (see equivalent definitions in~\cite[Corollary~6.6]{DOG}):
$$
\hspace*{12.5mm}\begin{array}{ll}
E_G(g)\! \!\! & =\{f\in G\,|\, \exists\,  n\in \mathbb{N}:  f^{-1}g^nf=g^{\pm n}\}\vspace*{3mm}\\
\! \!\! & =\{f\in G\,|\, \exists\,  k,m\in \mathbb{Z}\setminus \{0\}:  f^{-1}g^kf=g^{m}\}.
\end{array}
$$

\begin{lemma}\label{elem_index} {\rm (see~\cite[Lemma 6.8]{Osin_1})}
Suppose that $G$ is a group acting acylindrically on a hyperbolic space $S$. Then there exists $L\in \mathbb{N}$
such that for every loxodromic element $g\in G$, $E_G(g)$ contains a normal infinite cyclic subgroup
of index~$L$.
\end{lemma}



\begin{definition}\label{Def_L}
{\rm Let $G$ be a group and $X$ be a generating set of $G$.
For any two elements $u,v\in G$, we choose a geodesic path $[u,v]$ in $\Gamma(G,X)$ from $u$ to~$v$ so that
$w[u,v]=[wu,wv]$ for any $w\in G$.
With any element $x\in G$ and any loxodromic element $g\in G$, we associate the bi-infinite quasi-geodesic
$$
L(x,g)=\overset{\infty}{\underset{i=-\infty}{\cup}}x[g^i,g^{i+1}].
$$
We have $L(x,g)=x\, L(1,g)$. The path $L(1,g)$ is called the {\it quasi-geodesic associated with $g$}.
}
\end{definition}

\begin{corollary}\label{qg} {\rm (\cite[Corollary 2.12]{Bog_0})} Let $G$ be a group and $X$ be a generating set of $G$. Suppose that the Cayley graph $\Gamma(G,X)$ is hyperbolic and acylindrical.
Then there exist $\varkappa\geqslant 1$ and $\varepsilon\geqslant 0$ such that the following holds:

If an element $g\in G$ is loxodromic and shortest in its conjugacy class, then the quasi-geodesic $L(1,g)$
associated with $g$ is a $(\varkappa,\varepsilon)$-quasi-geodesic.
\end{corollary}

We will use the following technical lemmas from~\cite{Bog_1}.

\begin{lemma}\label{quasi-geod_1} {\rm (see~\cite[Lemma 4.7]{Bog_1})}
Let $G$ be a group and $X$ be a generating set of $G$. Suppose that the Cayley graph $\Gamma(G,X)$ is
hyperbolic and acylindrical.
Then there exist real numbers $\varkappa\geqslant 1, \varepsilon_0\geqslant 0$
and a number $n_0\in \mathbb{N}$ with the following property.

Suppose that $n\geqslant n_0$ and $c\in G$ is a loxodromic element.
Let $S(c)$ be the set of shortest elements in the conjugacy class
of $c$ and let $g\in G$ be a shortest element for which there exists $c_1\in S(c)$ with $c=g^{-1}c_1g$.
Then any path $p_0p_1\dots p_np_{n+1}$ in $\Gamma(G,X)$, where $p_0, p_1,\dots,p_n,p_{n+1}$
are geodesics with labels representing $g^{-1},c_1,\dots, c_1,g$, is a $(\varkappa, \varepsilon_0)$-quasi-geodesic.
In particular,
$$
|c^n|_X\geqslant \frac{1}{\varkappa}\bigl(n|c_1|_X+2|g|_X\bigr)-\varepsilon_0.
$$
\end{lemma}

\subsection{Stable norm}

Let $G$ be a group and $X$ is a generating set of $G$.
Recall that the {\it stable norm} of an element $g\in G$ with respect to a generating set $X$ is defined as $$||g||_X=\underset{n\rightarrow \infty}{\lim}\frac{|g^n|_X}{n},$$
see~\cite{CDP}. It is easy to check that this number is well-defined, that it is a conjugacy invariant, and that $||g^k||_X=|k|\cdot||g||_X$ for all $k\in \mathbb{Z}$.

Bowditch~\cite[Lemma 2.2]{Bowditch} proved that every element of a group acting acylindrically on a hyperbolic space is either elliptic or loxodromic (see a more general statement in~\cite[Theorem 1.1]{Osin_1}).
Moreover, he proved there
that the infimum of the set of stable norms of all loxodromic elements for such an action is larger than zero
(we assume that $\inf \emptyset=+\infty$).

\begin{lemma}\label{norm_length} Let $G$ be a group and $X$ be a generating set of $G$. Suppose that the Cayley graph $\Gamma(G,X)$ is hyperbolic and acylindrical.
For any loxodromic element $a\in G$, which is shortest in its conjugacy class, we have
$$
||a||_X\geqslant \frac{|a|_X}{\varkappa},\eqno{(2.1)}
$$
where $\varkappa\geqslant 1$ is the universal constant from Corollary~\ref{qg}.
\end{lemma}

{\it Proof.} By Corollary~\ref{qg}, there exist universal constants $\varkappa\geqslant 1$ and $\varepsilon\geqslant 0$ such that the path $L(1,a)$ is a $(\varkappa,\varepsilon)$-quasi-geodesic. Then,
for any natural $n$, we have
$$
|a^n|_X\geqslant \frac{\ell(a^n)-\varepsilon}{\varkappa}=\frac{n|a|_X-\varepsilon}{\varkappa}.
$$
Therefore
$$
||a||_X=\underset{n\rightarrow \infty}{\lim}\frac{|a^n|_X}{n}\geqslant \frac{|a|_X}{\varkappa}.
$$
\hfill $\Box$

\section{Proof of Theorem~\ref{third_main theorem}}

Theorem~\ref{third_main theorem} will be deduced from the following two lemmas, which say (simplified) that
an element $h\in G$ can be conjugate to a shortest representative by a element $g$, whose length
is bounded by a linear function of the length of $h$.
The first lemma (about loxodromic $h$) follows directly from Lemma~\ref{quasi-geod_1}, while the second one
(about elliptic $h$) seems to be not evident and needs an extended proof.

\begin{lemma}\label{loxo}
Let $G$ be an acylindrically hyperbolic group with respect to a generating set $X$. Then for any $h\in {\rm Lox}(G,X)$, there exists $g\in G$ such that $g h g^{-1}$ is a shortest element in the conjugacy class of $h$ and $|g|_X\leqslant K|h|_X$, where $K>0$ is a universal constant depending on the acylindricity data of the pair $(G,X)$.
\end{lemma}


{\it Proof.} By Lemma~\ref{quasi-geod_1}, there exists universal constants $n_0$, $\varkappa$ and $\varepsilon_0$ such that
$|g|_X\leqslant \frac{1}{2}\varkappa (|h^{n_0}|_X+\varepsilon_0)$. Then the statement holds for $K=\frac{1}{2}\varkappa n_0+\varepsilon_0$.\hfill $\Box$

\begin{lemma}\label{elli}
Let $G$ be an acylindrically hyperbolic group with respect to a generating set $X$. Then for any $h\in {\rm Ell}(G,X)$, there exists $g\in G$ such that $g\langle h\rangle g^{-1}\subseteq  {\rm \bf B}_1(8\delta+1)$ and $|g|_X\leqslant K|h|_X$, where $K>0$ is a universal constant depending on the acylindricity data of the pair $(G,X)$.
\end{lemma}


{\it Proof.} It is known that there exists $g\in G$ such that $g\langle h\rangle g^{-1}\subseteq  {\rm \bf B}_1(4\delta+1)$,
see~\cite[Corollary 6.7]{Osin_1}  (the proof there utilizes the proof of~\cite[Part III $\Gamma$, Theorem 3.2]{BH}). We start with some $g$ satisfying this property and modify it to get a (possibly) other $g$ with the desired length. For any integer $i$ we denote $h_i=gh^ig^{-1}$. Clearly $h_i=h_1^i$
and
$$
|h_i|_X\leqslant 4\delta+1.\eqno{(3.1)}
$$
We choose a geodesic path $[A,B]$ in $\Gamma(G,X)$ from $A=1$ to $B=g$.
For any $i\in \mathbb{Z}$ we consider the geodesic path $[A_i,B_i]=h_i[A,B]$
and choose geodesic paths $[A,A_i]$ and $[B,B_i]$. Note that
the paths $[A,B]$ and $[A_i,B_i]$ are both labeled by $g$ and the paths
$[A,A_i]$ and $[B,B_i]$ are labeled by $h_i$ and $h^i$, respectively, see Fig. 1.

\vspace*{-30mm}
\hspace*{0mm}
\includegraphics[scale=0.8]{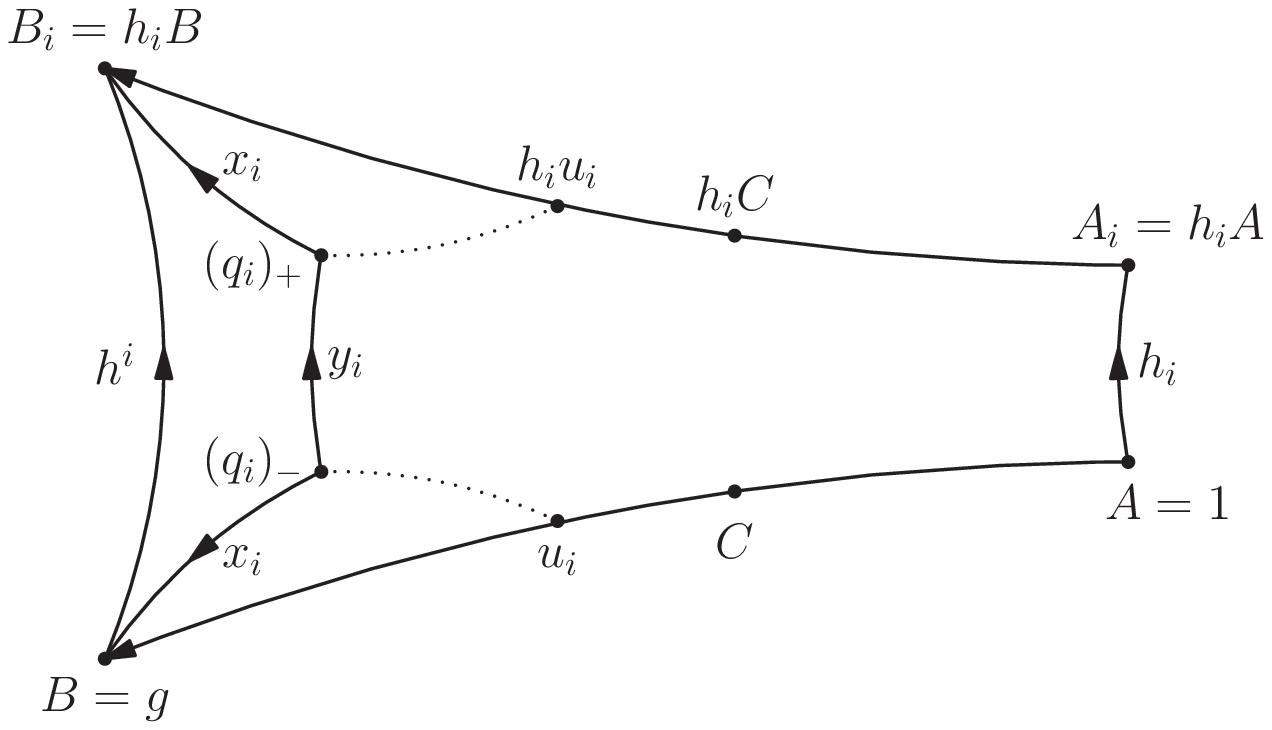}

\vspace*{-135mm}

\begin{center}
Fig. 1. Illustration to the proof of the main theorem.
\end{center}

If $|g|_X< R(8\delta+3)$, we are done with $K=R(8\delta+3)$. Therefore we assume that
$$
d(A,B)=|g|_X\geqslant R(8\delta+3).\eqno{(3.2)}
$$
Let
$$
I=\{i\in \mathbb{Z}\,|\, |h^i|_X\leqslant 8\delta+3\}.\eqno{(3.3)}
$$

{\bf Claim 1.} We have
$
\#\{h^i\,|\, i\in I\}\leqslant N(8\delta+3).
$

\medskip

{\it Proof.}
For any $i\in \mathbb{Z}$, we have
$$
d(A,h_iA)=|h_i|_X\overset{(3.1)}{\leqslant}4\delta+1
$$
and for any $i\in I$ we have
$$
d(B,h_iB)=d(g,h_ig)=|g^{-1}h_ig|_X=|h^i|_X\overset{(3.3)}{\leqslant} 8\delta+3.
$$
From this, (3.2) and the definition of the acylindrical action we obtain the statement.\hfill $\Box$

\medskip

Now consider $i\in I^c$, where $I^c=\mathbb{Z}\setminus I$.
By Lemma~\ref{quasi-geod_2} applied to $a=h^i$ and $b=h_i$, there exist
$x_i,y_i\in G$ such that
$$
h^i=x_i^{-1}y_ix_i,
$$
$$
|y_i|_X\leqslant 8\delta+3,\eqno{(3.4)}
$$
and any path $\ell_i=p_iq_ir_i$ in $\Gamma(G,X)$, where $p_i,q_i,r_i$ are geodesics
with labels $x_i^{-1}$, $y_i$, $x_i$ is a $(1,\epsilon_1)$-quasi-geodesic. Here $\epsilon_1=\epsilon_1(\delta)$
is a universal constant. In particular, we have
$$
2|x_i|_X+|y_i|_X\leqslant |h^i|_X+\epsilon_1,\eqno{(3.5)}
$$

Since the labels of $\ell_i$ and $[B,B_i]$ are both equal to $h^i$, we can choose $\ell_i$ so that $(\ell_i)_{-}=B$ and $(\ell_i)_{+}=B_i$, see Fig. 1.
Observe that
$$
h_i\cdot(q_i)_{-}=gh^ig^{-1}\cdot gx_i^{-1}=(q_i)_{+}.\eqno{(3.6)}
$$
Since the label of $q_i$ is $y_i$, we deduce from (3.4) that
$$
d((q_i)_{-},(q_i)_{+})\leqslant 8\delta+3.\eqno{(3.7)}
$$

\medskip

{\bf Claim 2.} There exists a constant $\epsilon_2>0$ depending only on $\delta$ such that the following holds.
For any $i\in I^c$, there exists a point $u_i\in [A,B]$ such that
$$
d((q_i)_{-}, u_i)\leqslant \epsilon_2.\eqno{(3.8)}
$$
and
$$
d(u_i,h_iu_i)\leqslant 2\epsilon_2+8\delta+3.\eqno{(3.9)}
$$

\medskip

{\it Proof.} We set  $\mu_1=\mu(\delta,1,\epsilon_1)$, where the function $\mu$ is defined in  Lemma~\ref{Hausdorff}
and prove that the statement is valid for
$
\epsilon_2=\mu_1+10\delta+3.
$

We prove the first statement. Recall that $\ell_i=p_iq_ir_i$ is a $(1,\epsilon_1)$-quasi-geodesic with endpoints $B,B_i$. Then, by Lemma~\ref{Hausdorff}, there exists a point $w_i\in [B,B_i]$ such that $d(w_i,(q_i)_{-})\leqslant \mu_1$. By Lemma~\ref{polygon}, $w_i$ is at distance at most $2\delta$
from the union of three sides $[B,A]$, $[A,A_i]$, $[A_i,B_i]$.

{\it Case 1.} Suppose that there exists $z_i\in [B,A]$ such that $d(w_i,z_i)\leqslant 2\delta$.
Then
$$
d((q_i)_{-}, z_i)\leqslant d((q_i)_{-},w_i)+d(w_i,z_i)\leqslant \mu_1+2\delta<\epsilon_2,
$$
and we are done with $u_i=z_i$.

\medskip

{\it Case 2.} Suppose that there exists $z_i\in [A,A_i]$ such that $d(w_i,z_i)\leqslant 2\delta$.
Then
$$d((q_i)_{-}, A)\leqslant d((q_i)_{-},w_i)+d(w_i,z_i)+d(z_i,A)\leqslant \mu_1 +2\delta+(4\delta+1)<\epsilon_2,$$ and we are done with $u_i=A$.

\medskip

{\it Case 3.} Suppose that there exists $z_i\in [A_i,B_i]$ such that $d(w_i,z_i)\leqslant 2\delta$.
We set $u_i=h_i^{-1}z_i$. Then $u_i\in [A,B]$ and we have
$$
\begin{array}{ll}
d((q_i)_{-}, u_i) & \overset{(3.6)}{=}d(h_i^{-1}(q_i)_{+},h_i^{-1}z_i)\vspace*{2mm}\\
                  & =d((q_i)_{+},z_i)\vspace*{2mm}\\
& \leqslant d((q_i)_{+},(q_i)_{-})+d((q_i)_{-},w_i)+d(w_i,z_i)\vspace*{2mm}\\
& \overset{(3.7)}{\leqslant} (8\delta+3)+\mu_1+2\delta=\epsilon_2.
\end{array}
$$
This completes the proof of the first statement. Now we prove the second statement:
$$
\begin{array}{ll}
d(u_i,h_iu_i)&\leqslant d(u_i,(q_i)_{-})+d((q_i)_{-},(q_i)_{+})+d((q_i)_{+},h_iu_i)\vspace*{2mm}\\
& \overset{(3.5)}{=} d(u_i,(q_i)_{-})+d((q_i)_{-},(q_i)_{+})+d(h_i(q_i)_{-},h_iu_i)\vspace*{2mm}\\
& \underset{(3.7)}{\overset{(3.6)}{\leqslant}} \epsilon_2+(8\delta+3)+\epsilon_2.\vspace*{2mm}\\
\end{array}
$$

\hfill $\Box$

Now we define the set
$$
J =\{i\in I^c\,|\, d(A,u_i)>R(8\delta+1)+2\epsilon_2+16\delta+5\}.
$$

{\bf Claim 3}. We have
$
\#\{h^i\,|\, i\in J\}\leqslant N(8\delta+1).
$

\medskip

{\it Proof.} We assume that $J\neq \varnothing$. Then there exists a point $C\in [A,B]$ such that
$$
d(A,C)=R(8\delta+1)+6\delta+2,\eqno{(3.10)}
$$
and we have $C\in [A,u_i]$ for any $i\in J$. First we prove that
$$
d(C,h_iC)\leqslant 8\delta+1.\eqno{(3.11)}
$$
For that we apply Lemma~\ref{equidist} to the geodesic paths $[A,u_i]$ and $[h_iA,h_iu_i]$ and the points
$C$ and $h_iC$. The assumptions of this lemma are satisfied:

\begin{enumerate}
\item[(a)] $d(C,A)\overset{(3.10)}{>} 6\delta+2\overset{(3.1)}{\geqslant} d(A,A_i)+2\delta$.

\item[(b)] $d(C,u_i)=d(A,u_i)-d(A,C)\overset{(3.10)}{>}2\epsilon_2+10\delta+3\overset{(3.9)}{\geqslant} d(u_i,h_iu_i)+2\delta$.

\item[(c)] $d(A,C)=d(h_iA,h_iC)$.
\end{enumerate}

By this lemma, $d(C,h_iC)\leqslant 4\delta+d(A,h_iA)\leqslant 4\delta+(4\delta+1)$ that proves (3.11).
By (3.1) we have $d(A,h_iA)\leqslant 4\delta+1$ and by (3.10) we have $d(C,A)>R(8\delta+1)$.
From this, (3.11) and the definition of the acylindrical action we obtain the statement.
\hfill $\Box$

\medskip

Now we are ready to complete the proof of the statement. It follows from Claims~1 and~3 that
$$
\#\{h^i\,|\, i\in I\cup J \}\leqslant n,\hspace*{3mm}{\rm where}\hspace*{3
mm}n=N(8\delta+3)+N(8\delta+1)\eqno{(3.12)}
$$

{\it Case 1.} Suppose  $\#\langle h\rangle\leqslant n$.

Let $\mathcal{M}=\max\{|h^i|_X: 1\leqslant i\leqslant n\}$. Note that $\mathcal{M}\leqslant n|h|_X$.
If $|g|_X\leqslant \mathcal{M}+8\delta+2$, we are done.
Suppose that $|g|_X> \mathcal{M}+8\delta+2$. Let $C$ be the point on the side $[A,B]$ such that $$
d(C,B)=\mathcal{M}+2\delta+1.\eqno{(3.13)}
$$
Then $d(C,A)> 6\delta+1$.
It follows that the distance from $C$ to $[A,A_i]\cup [B,B_i]$ is larger than $2\delta$.
We set $C_i=h_iC$. Then, by Lemma~\ref{equidist},
$$
d(C,C_i)\leqslant 8\delta+1.\eqno{(3.14)}
$$

Let $g_1$ be the label of the path $[C,B]$ (and hence of the path $[C_i,B_i]$).
The concatenation of the paths $[C,B]$, $[B,B_i]$, $[B_i,C_i]$ has the same endpoints as the geodesic path $[C,C_i]$. Therefore the label (in $G$) of the path $[C,C_i]$ is $g_1h^ig_1^{-1}$.
Using (3.14), we obtain $|g_1h^ig_1^{-1}|_X\leqslant 8\delta+1$ for any $i$.
Using (3.12) and (3.13), we deduce
$$
\begin{array}{ll}
|g_1|_X=d(C,B) & =\mathcal{M}+2\delta+1\vspace*{2mm}\\
& \leqslant n|h|_X+2\delta+1\vspace*{2mm}\\
& \leqslant (N(8\delta+3)+N(8\delta+1))|h|_X+2\delta+1.
\end{array}
$$
This completes the proof in this case.

\medskip

{\it Case 2.} Suppose  $\#\langle h\rangle> n$.

By (3.12), one of the elements $1,h,h^2,\dots ,h^n$ does not lie in the set $\{h^i\,|\, i\in I\cup J \}$.
Then there exists $0\leqslant i\leqslant n$ such that $i\in I^c\setminus J$.
In particular, $d((q_i)_{-},u_i)\leqslant \epsilon_2$ and $d(A,u_i)\leqslant R(8\delta+1)+2\epsilon_2+16\delta+5$. Then
$$
\begin{array}{ll}
|g|_X & \leqslant d(B,(q_i)_{-})+d((q_i)_{-},u_i)+d(u_i,A)\vspace*{2mm}\\
& \leqslant |x_i|_X+\epsilon_2+(R(8\delta+1)+2\epsilon_2+16\delta+5).
\end{array}
$$
Finally we note that
$$
|x_i|_X\overset{(3.5)}{\leqslant} \frac{1}{2}(|h^i|_X+\epsilon_1)\leqslant \frac{1}{2}(n|h|_X+\epsilon_1).
$$
This completes the proof.\hfill $\Box$

\medskip

{\it Proof of Theorem~\ref{third_main theorem}.}
(1) Suppose that $h_1,h_2$ are loxodromic elements. By Lemma~\ref{loxo}, we may reduce the proof to the case that $h_1,h_2$ are shortest in their conjugacy class.
Let $g\in G$ be an arbitrary element such that $h_1=gh_2g^{-1}$.
We make two observations about the quasi-geodesics $L(1,h_1)$ and $L(g,h_2)$.

\begin{enumerate}
\item[(a)] Since $h_1=gh_2g^{-1}$, the Hausdorff distance between $L(1,h_1)$ and $L(g,h_2)$ is at most $|g|_X+\max\{|h_1|_X, |h_2|_X\}$. Therefore the limit points of these quasi-geodesics coincide.

\item[(b)] Since $h_1$ and $h_2$ are shortest in their conjugacy class, both $L(1,h_1)$ and $L(g,h_2)$ are $(\varkappa,\varepsilon)$-quasi-geodesics, where $\varkappa$ and $\varepsilon$ are
universal constants from Corollary~\ref{qg}.
\end{enumerate}

It follows from (a) and (b) that the Hausdorff distance between $L(1,h_1)$ and $L(g,h_2)$ is at most $k=\mu(\delta,\varkappa,\varepsilon$), see Lemma~\ref{HausdorffQuasiGeod}.
In particular, there exists a point $z\in L(g,h_2)$ such that
$d(1,z)\leqslant k$. Let $t=gh_2^i$ be the phase point on $L(g,h_2)$ which is nearest to $z$.
In particular, $d(t,z)\leqslant |h_2|_X$. Then
$$
|t|_X=d(1,t)\leqslant d(1,z)+d(z,t)\leqslant k+|h_2|_X.
$$
Moreover, $h_1=th_2t^{-1}$, and we are done.

\medskip

(2) Suppose that $h_1,h_2$ are elliptic elements.
By Lemma~\ref{elli}, we may reduce the proof to the case that the subgroups $\langle h_i\rangle$, $i=1,2$, lie in the ball $B_1(8\delta+1)$. Let $g\in G$ be an element such that $h_1=gh_2g^{-1}$.
Since the orders of $h_i$ are larger than $N(8\delta+1)$ (by assumption), it follows from the definition of acylindricity that $|g|_X\leqslant R(8\delta+1)$.
\hfill $\Box$

\section{An extension of the periodicity theorem from~\cite{Bog_0}}

The main result of this section is Theorem~\ref{acylindric_1}, which is used in Section~6.
It slightly extends the periodicity theorem from~\cite{Bog_0}, see Theorem~\ref{acylindric} below.
Both theorems can be easily formulated in the case of free groups:

Let $a,b$ be two cyclically reduced words in the free group $F$ with basis $X$.
If the bi-infinite words $L(a)=\dots aaa\dots $ and $L(b)=\dots bbb\dots$ have a common subword of length $|a|+|b|$,
then some cyclic permutations of $a$ and $b$ are positive powers of some word $c$.


For the case of acylindrically hyperbolic groups, we recall some notions.
Suppose that the Cayley graph $\Gamma(G,X)$
is hyperbolic and that $G$ acts acylindrically on $\Gamma(G,X)$. In~\cite[Lemma 2.2]{Bowditch}
Bowditch proved that the infimum of stable norms (see Section~2) of all loxodromic
elements of $G$ with respect to $X$ is a positive number.
We denote this number by ${\bf inj}(G,X)$ and call it the {\it injectivity radius} of $G$ with respect to $X$.







\begin{definition}\label{L}
{\rm Let $G$ be a group and $X$ a generating set of $G$. The right Cayley graph of $G$
with respect to $X$ is denoted by $\Gamma(G,X)$. For any two elements $u,v\in G$, we choose a geodesic path $[u,v]$ in $\Gamma(G,X)$ from $u$ to $v$ so that
$w[u,v]=[wu,wv]$ for any $w\in G$.
With any element $x\in G$ and any element $g\in G$ of infinite order, we associate the bi-infinite path
$L(x,g)=\dots p_{-1}p_0p_1\dots$, where $p_n=[xg^n,xg^{n+1}]$, $n\in \mathbb{Z}$.
The paths $p_n$ are called {\it $g$-periods} of $L(x,g)$. For a subpath $p\subset L(x,g)$  and a number $k\in \mathbb{N}$, we say that the path $p$ {\it contains $k$ $g$-periods} if there exists $n\in \mathbb{Z}$ such that $p_np_{n+1}\dots p_{n+k-1}$ is a subpath of $p$.
The vertices $xg^n$, $n\in \mathbb{Z}$, are called the {\it phase} vertices of $L(x,g)$.
Note that $L(x,g)=x\, L(1,g)$.
}
\end{definition}


\begin{theorem}{\rm (see~\cite[Theorem 1.4]{Bog_0})}\label{acylindric}
Let $G$ be a group and $X$ a generating set of $G$. Suppose that the Cayley graph $\Gamma(G,X)$
is hyperbolic and that $G$ acts acylindrically on $\Gamma(G,X)$.
Then there exists a constant $\mathcal{C}>0$ such that the following holds.

Let $a,b\in G$ be two loxodromic elements which are shortest in their conjugacy classes
and such that $|a|_X\geqslant |b|_X$.
Let $x,y\in G$ be arbitrary elements and $r$ an arbitrary non-negative real number. We set $f(r)=\frac{2r}{{\bf inj}(G,X)}+\mathcal{C}$.


Suppose that $p\subset L(x,a)$ and $q\subset L(y,b)$ are subpaths such that $d(p_{-},q_{-})\leqslant r$, $d(p_{+},q_{+})\leqslant r$,
and $p$ contains at least $f(r)$ $a$-periods.
Then there exist nonzero integers $s,t$ such that
$$
(y^{-1}x)a^s(x^{-1}y)=b^t.
$$


\end{theorem}

\vspace*{-20mm}
\hspace*{25.5mm}
\includegraphics[scale=0.55]{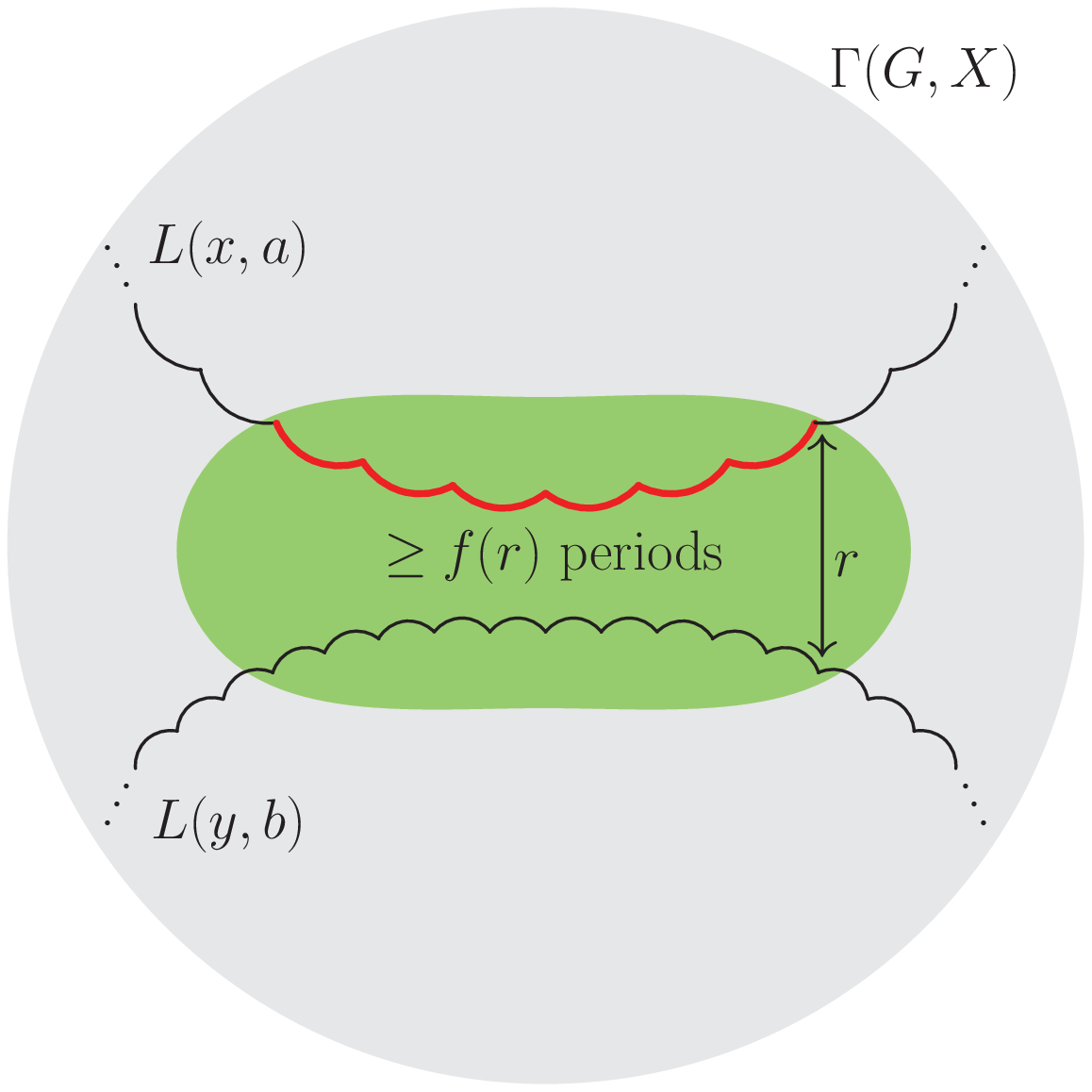}

\vspace*{-72mm}

\begin{center}
Fig. 2. Illustration to Theorem~\ref{acylindric}.
\end{center}

The following theorem says that taking an appropriate linear function $F$ instead of $f$, we can guarantee that
both numbers $s$ and $t$ are positive.

\begin{theorem}\label{acylindric_1}
Let $G$ be a group and $X$ a generating set of $G$. Suppose that the Cayley graph $\Gamma(G,X)$
is hyperbolic and that $G$ acts acylindrically on $\Gamma(G,X)$.
Then there exists a a linear function $F:\mathbb{R}\rightarrow \mathbb{R}$ with constants depending only on $(G,X)$ such that
the following holds.

Let $a,b\in G$ be two loxodromic elements which are shortest in their conjugacy classes
and such that $|a|_X\geqslant |b|_X$.
Let $x,y\in G$ be arbitrary elements and $r$ an arbitrary non-negative real number.
Suppose that $p\subset L(x,a)$ and $q\subset L(y,b)$ are subpaths such that $d(p_{-},q_{-})\leqslant r$, $d(p_{+},q_{+})\leqslant r$,
and $p$ contains at least $F(r)$\break $a$-periods.
Then there exist positive integers $s,t$ such that
$$
(y^{-1}x)a^s(x^{-1}y)=b^t.
$$
\end{theorem}


{\it Proof.}
We set $F(r)= f(r)+\varkappa (4r+4\delta+5\mu)+\varepsilon+1$, where $\delta$ is the hyperbolicity constant of the Cayley graph $\Gamma(G,X)$,  $\varkappa$ and $\varepsilon$ are from Corollary~\ref{qg}, and $\mu=\mu(\delta,\varkappa,\varepsilon)$ is from Lemma~\ref{Hausdorff}.
For brevity, we set $\mathcal{F}=\lfloor F(r)\rfloor $.

First we show that, without loss of generality, we may assume that $p_{-}$ and $q_{-}$ are phase vertices of $L(x,a)$
and $L(y,b)$, respectively.

Let $A$ and $B$ be the leftmost phase vertices of $p$ and $q$, respectively.
We set $a_1=cac^{-1}$, where $c$ is the subpath of the $a$-period from $p_{-}$ to $A$ and we set $b_1=dbd^{-1}$, where $d$ is the subpath of the $b$-period from $q_{-}$ to $B$, see Figure 3.

\vspace*{-30mm}
\hspace*{-0mm}
\includegraphics[scale=0.8]{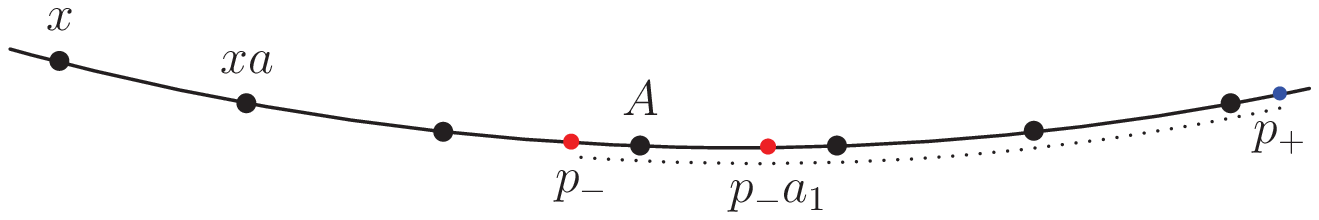}

\vspace*{-185mm}

\begin{center}
Fig. 3. Reduction $L(x,a)=L(p_{-},a_1)$.
\end{center}

Then $|a_1|_X=|a|_X$ and $L(x,a)=L(p_{-},a_1)$ and also $|b_1|_X=|b|_X$ and $L(y,b)=L(q_{-},b_1)$.
Note that $p_{-}$ is a phase vertex of $L(p_{-},a_1)$ and $q_{-}$ is a phase vertex of $L(q_{-},b_1)$.
Suppose we have proved that there exist positive integers $s,t$ such that
$$
(q_{-}^{-1}p_{-})a_1^s(p_{-}^{-1}q_{-})=b_1^t.
$$
Substituting $a_1=cac^{-1}$, $b_1=dbd^{-1}$, $A=p_{-}c$ and $B=q_{-}d$, we deduce
$$
(B^{-1}A)a^s(A^{-1}B)=b^t.
$$
Since $A$ is a phase vertex of $L(x,a)$, we have $A=xa^i$ for some $i\in \mathbb{Z}$.
Analogously we have $B=yb^j$ for some $j\in \mathbb{Z}$. This implies $(y^{-1}x)a^s(x^{-1}y)=b^t.$

\medskip

Thus, without loss of generality, we assume that $p_{-}$ and $q_{-}$ are phase vertices of $L(x,a)$ and $L(y,b)$, respectively. Then $L(x,a)=L(p_{-},a)$ and $L(y,b)=L(q_{-},b)$, and by Theorem~\ref{acylindric} we have
$$
(q_{-}^{-1}p_{-})a^s(p_{-}^{-1}q_{-})=b^t\eqno{(4.1)}
$$
for some nonzero integers $s,t$. We may assume that $s>0$. Suppose that $t<0$.
We set $C=p_{-}a^{s\mathcal{F}}$ and $D=q_{-}b^{t\mathcal{F}}$. Then $C$ lies on $L(p_{-},a)$ to the right from $p_{+}$
and $D$ lies on $L(q_{-},b)$ to the left from $q_{-}$, see Figure 4.

\vspace*{-30mm}
\hspace*{-1mm}
\includegraphics[scale=0.8]{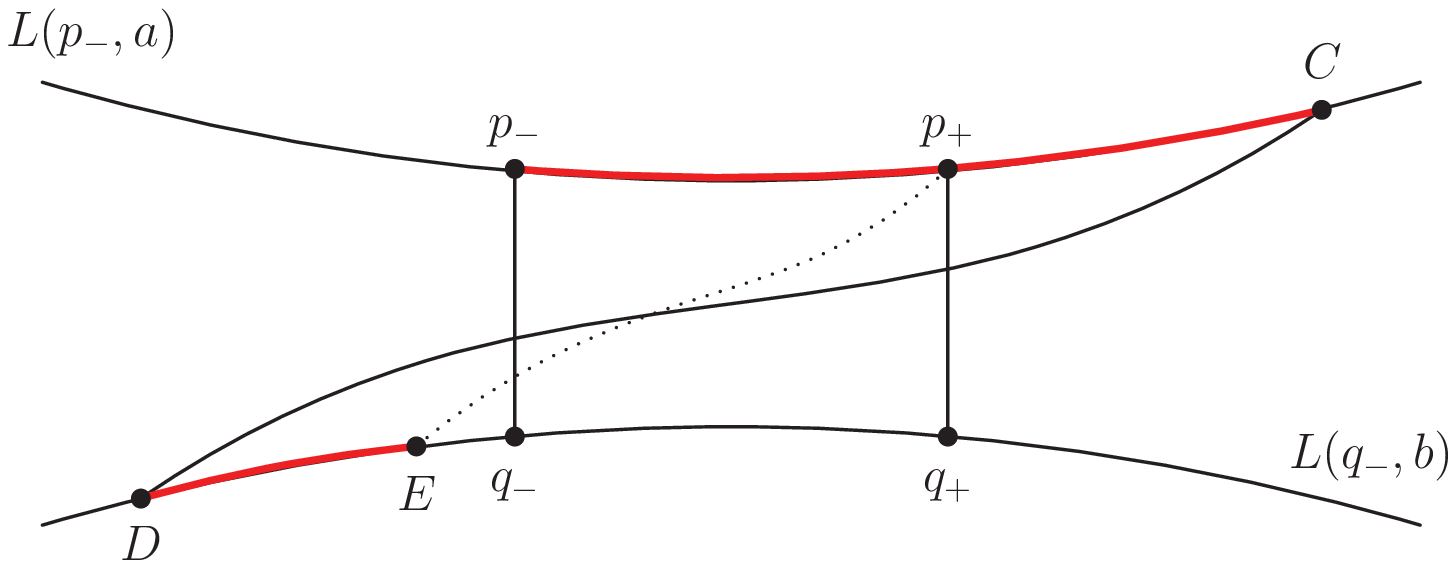}

\vspace*{-155mm}

\begin{center}
Fig. 4. The case $s>0$ and $t<0$.
\end{center}

\noindent
Let $u$ be the subpath of $L(p_{-},a)$ from $p_{-}$ to $C$ and let
$v$ be the subpath of $L(q_{-},b^{-1})$ from $q_{-}$ to $D$.
We have $d(u_{-},v_{-})=d(p_{-},q_{-})\leqslant r$ by assumption in the theorem and we have $d(u_{+},v_{+})\leqslant r$ since
$$
d(u_{+},v_{+})=d(p_{-}a^{s\mathcal{F}}, q_{-}b^{t\mathcal{F}})=d(1,a^{-s\mathcal{F}}p_{-}^{-1}q_{-}b^{t\mathcal{F}})\overset{(4.1)}{=}d(1, p_{-}^{-1}q_{-})=d(p_{-},q_{-})\leqslant r.
$$
By Corollary~\ref{close_qg}, there exists a point $E\in v$ such that
$$
d(p_{+},E)\leqslant r+2(\delta+\mu).\eqno{(4.2)}
$$
We have
$$
d(E,q_{-})\geqslant d(p_{-},p_{+})-d(p_{-},q_{-})-d(p_{+},E)\geqslant d(p_{-},p_{+})-r-(r+2(\delta+\mu)).\eqno{(4.3)}
$$
The point $q_{-}$ lies on the $(\varkappa,\varepsilon)$-quasi-geodesic $L(y,b)$ between the points $E$ and $q_{+}$. Therefore, by Lemma~\ref{Hausdorff}, there
exists a point $q_{-}'\in [E,q_{+}]$ such that $d(q_{-},q_{-}')\leqslant \mu$.
Then
$$
\begin{array}{ll}
d(E,q_{-}) & \leqslant d(E,q_{-}')+d(q_{-}', q_{-})\vspace*{2mm}\\
& \leqslant d(E,q_{+})+\mu\vspace*{2mm}\\
& \leqslant d(E,p_{+})+d(p_{+},q_{+})+\mu\vspace*{2mm}\\
& \overset{(4.2)}{\leqslant} (r+2(\delta+\mu))+r+\mu.
\end{array}\eqno{(4.4)}
$$

It follows from (4.3) and (4.4) that

$$
d(p_{-},p_{+})\leqslant 4r+4\delta+5\mu.
$$
On the other hand,
$$
\mathcal{F}\leqslant \mathcal{F}|a|_X\leqslant \ell(p)\leqslant \varkappa d(p_{-},p_{+})+\varepsilon\leqslant
\varkappa (4r+4\delta+5\mu)+\varepsilon
$$
that contradicts the definition of $\mathcal{F}$ at the beginning of the proof.
Thus the assumption $t<0$ is not valid.\hfill $\Box$

\medskip

\noindent
{\bf Notation.}
For any subpath $p\subset L(x,g)$ let $N(p)$ the number of $g$-periods containing in $p$.
In Section 6, we will use the following easy observation.
$$
N(p)|g|_X+2|g|_X\geqslant \ell(p)\geqslant d(p_{-},p_{+})\geqslant |g^{N(p)}|_X-2|g|_X.\eqno{(4.5)}
$$

\section{Indices}



We need the following generalization of the notion the least common multiple of two nonzero integers.
In the case of $\mathbb{Z}$ the index introduced in the following definition coincides
with the index of the ideal $(a,b)$ in the ideal $(a)$.

\begin{definition}
Let $G$ be a group and let $[a], [b]$ be two conjugacy classes of elements $a,b\in G$ of infinite order.
Suppose that $a,b$ are commensurable. Then, by definition, there exist nonzero integers $k,\ell$ such that the conjugacy classes of $a^k$ and  $b^{\ell}$ coincide. We take minimal $k>0$ with this property and call the conjugacy class $of a^k$ the {\it least common multiple of the conjugacy classes} of $a$ and $b$, and we denote it by $[a]\vee [b]$. The number $k$ is called the {\it index of $[a]\vee [b]$ with respect to $[a]$} and is denoted by ${\bf Ind}_{[a]}([a]\vee [b])$.  Thus,
$$
{\bf Ind}_{[a]}([a]\vee [b]):=\min\{ k>0\,|\, \exists\, s: a^k\sim b^s \}.
$$
\end{definition}

\begin{remark} The conjugacy class $[a]\vee [b]$ does not depend of the choice of $a$ and $b$ in their conjugacy classes. The following lemma implies that if $a$ and $b$ are loxodromic elements of an acylindrically hyperbolic group $G$, then $[a]\vee [b]=\pm ([b]\vee [a])$. It also gives an estimation of ${\bf Ind}_{[a]}([a]\vee [b])$ via the stable norm of $b$.
\end{remark}

In the following lemmas $L$ is the constant from Lemma~\ref{elem_index}.

\begin{lemma}\label{estim_k} Let $G$ be an acylindrically hyperbolic group with respect to a generating system $X$. Let $a,b$ be two commensurable loxodromic elements of $G$. Denoting
$k={\bf Ind}_{[a]}([a]\vee [b])$ and $\ell={\bf Ind}_{[b]}([b]\vee [a])$, we have
$$
a^k\sim b^{\pm \ell},\eqno{(5.1)}
$$
$$
k\cdot ||a||_X=\ell\cdot ||b||_X,\eqno{(5.2)}
$$
$$
k\leqslant  \frac{L^2}{{\bf inj}(G,X)}\cdot ||b||_X.\eqno{(5.3)}
$$


\end{lemma}

\medskip

{\it Proof.} By definition we have $a^k\sim b^s$ and $b^{\ell}\sim a^t$ for some $s,t\in \mathbb{Z}$.
It follows $k\cdot||a||_X=|s|\cdot ||b||_X$ and $\ell\cdot ||b||_X=|t|\cdot||a||_X$.
Hence $k\ell=|s||t|$. By definition we have $k\leqslant |t|$ and $\ell\leqslant |s|$. This implies $s=\pm \ell$ and hence (5.1) and (5.2).

\medskip

We prove (5.3). By (5.1) we have $a^k=z^{-1}b^{\pm \ell}z$ for some $z\in G$.
It follows that $a\in E_G(z^{-1}bz)$.
Then, by Lemma~\ref{elem_index}, $a^L$ and $(z^{-1}bz)^L$ belong to the same infinite cyclic group. Let $c$ be a generator of this group.
Then $a^L=c^p$ and $(z^{-1}bz)^L=c^q$ for some nonzero integers $p,q$.
This implies
$$
a^{Lq}=z^{-1}b^{Lp}z.
$$
From this and the definition of $k$, we have $k\leqslant L|q|$. It remains to estimate $|q|$.\break
It follows from the definitions of stable norm and injectivity radius that
$$
L||b||_X=||b^L||_X=||c^q||_X=|q|\cdot ||c||_X\geqslant |q|\cdot {\bf inj}(G,X).
$$
Hence
$$
|q|\leqslant \frac{L||b||_X}{{\bf inj}(G,X)}.
$$
Substituting in the above established estimation $k\leqslant L|q|$, we complete the proof.
\hfill $\Box$.

The following lemma estimates possible nonzero exponents $s,t$ in the equation $z^{-1}a^sz=b^t$ with given $a,b,z\in G$, where $G$ is acylindrically hyperbolic and $a,b$ are loxodromic.

\begin{lemma}\label{nice_estimate}
Let $G$ be an acylindrically hyperbolic group with respect to a generating system $X$.
Let $a,b,z$ be elements of $G$, where $a$ and $b$ are loxodromic, such that $z^{-1}a^nz=b^m$ for some nonzero integers $n,m$.
Then we have $z^{-1}a^sz=b^t$ with the same $z$, where
$$
|s|= L\cdot {\bf Ind}_{[a]}([a]\vee [b])\hspace*{3mm}{\text{\rm and}}\hspace*{3mm} |t|= L\cdot {\bf Ind}_{[b]}([b]\vee [a]).
$$
Moreover, if $n,m$ are positive, then $s,t$ can be also chosen to be positive.
\end{lemma}

\medskip

{\it Proof.} We denote $k= {\bf Ind}_{[a]}([a]\vee [b])$ and $\ell={\bf Ind}_{[b]}([b]\vee [a])$.
By (5.1), there exists $z_1\in G$ such that
$$
z_1^{-1}a^{k}z_1=b^{\pm \ell}.
$$
From this and from the equation $z^{-1}a^nz=b^m$ we deduce
$$
z_1^{-1}a^{mk}z_1=b^{\pm m\ell}\hspace*{2mm}{\rm and}\hspace*{2mm}z^{-1}a^{n\ell}z=b^{m\ell}.
$$
We denote $e=zz_1^{-1}$. Then $ea^{mk}e^{-1}=a^{\pm n\ell}$, hence $e\in E_G(a)$.
By Lemma~\ref{elem_index}, $E_G(a)$ contains a normal infinite cyclic subgroup of index $L$.
It follows that
$e^{-1}a^Le=a^{\pm L}$.
Then
$$
z^{-1}a^{kL}z=z_1^{-1}e^{-1}a^{kL}ez_1=z_1^{-1}a^{\pm kL}z_1=b^{\pm \ell L}.
$$
This shows that the first statement is valid for $s=kL$ and $t=\pm \ell L$.

Now we prove the second statement. Suppose that both $n,m$ are positive.
From $z^{-1}a^nz=b^m$ and $z^{-1}a^sz=b^t$ follows $b^{ms}=b^{nt}$. Since $b$ has infinite order,
we have $ms=nt$, hence $s$ and $t$ have the same sign.
Changing the signs of $s$ and $t$ simultaneously, we may assume that both $s,t$ are positive. \hfill $\Box$



\section{An auxiliary lemma}

\begin{definition}
Let $G$ be a group and $g\in G$ be an element of infinite order. The set of elements of $G$ commensurable with $g$ is denoted by ${\text{\rm Com}}(g)$. Thus,
$$
{\text{\rm Com}}(g)=\{h\in G\,|\, g^t\hspace*{2mm} {\text {\rm is conjugate to}}\hspace*{2mm} h^s
\hspace*{2mm} {\text {\rm for some nonzero}}\hspace*{2mm} s,t\}.
$$
\end{definition}

\begin{lemma}\label{first_main lemma}
Let $G$ be an acylindrically hyperbolic group with respect to a generating set $X$.
Then there exists a constant $M>1$
such that for any exponential equation
$$
a_1g_1^{x_1}a_2g_2^{x_2}\dots a_ng_n^{x_n}=1\eqno{(6.1)}
$$
with constants $a_1,g_1,\dots,a_n,g_n$ from $G$ (where  $g_1,\dots, g_n$
are loxodromic and shortest in their conjugacy classes with respect to $X$) and variables $x_1,\dots,x_n$, if this equation is solvable over $\mathbb{Z}$, then there exists a solution $(k_1,\dots,k_n)$ with
$$
|k_j|\leqslant \Bigl(n^2+\overset{n}{\underset{i=1}\sum}\,\frac{|a_i|_X}{|g_j|_X}+
\underset{g_i\notin {\text{\rm Com}}(g_j)}\sum\,\frac{|g_i|_X}{|g_j|_X}+
\underset{g_i\in {\text{\rm Com}}(g_j)}{\sum}
{\bf Ind}_{[g_j]}([g_j]\vee [g_i])
\Bigr)\cdot M\eqno{(6.2)}
$$
for all $j=1,\dots,n$.

\end{lemma}

{\it Proof.}
Suppose that $(k_1,\dots,k_n)$ is a solution of equation (6.1) with minimal sum $|k_1|+\dots +|k_n|$. Because of symmetry, we estimate only $|k_1|$. In what follows we consider a polygon $\mathcal{P}$ in the Cayley graph $\Gamma(G,X)$ corresponding to the equation (6.1). More precisely, let $\mathcal{P}$ be a polygon in the Cayley graph $\Gamma(G,X)$ with consecutive sides
$p_1,q_1,p_2,q_2,\dots,p_n,q_n$ such that the sides $p_i$ are geodesics with the labels $a_i$
and the sides $q_i$ are quasi-geodesics consisting of $k_i$ consecutive geodesic paths labelled by $g_i$.
Note that each $q_i$ is a $(\varkappa,\varepsilon)$-quasi-geodesic path by Corollary~\ref{qg}.

By Lemmas~\ref{polygon} and~\ref{Hausdorff}, $q_1$ lies in the $\nu$-neighborhood of
the union of the other sides of $\mathcal{P}$, where
$$
\nu=(2n-2)\delta +2\mu.
$$
For $i=1,\dots,n$, let $p_i'$ be the maximal phase subpath of $q_1$ such that the endpoints of $p_i'$ are
at distance at most $\nu$ from $p_i$. Analogously, for  $i=2,\dots,n$,
let $q_i'$ be the maximal phase subpath of $q_1$ such that the endpoints of $q_i'$ are
at distance at most $\nu$ from $q_i$.

\vspace*{-30mm}
\hspace*{6mm}
\includegraphics[scale=0.8]{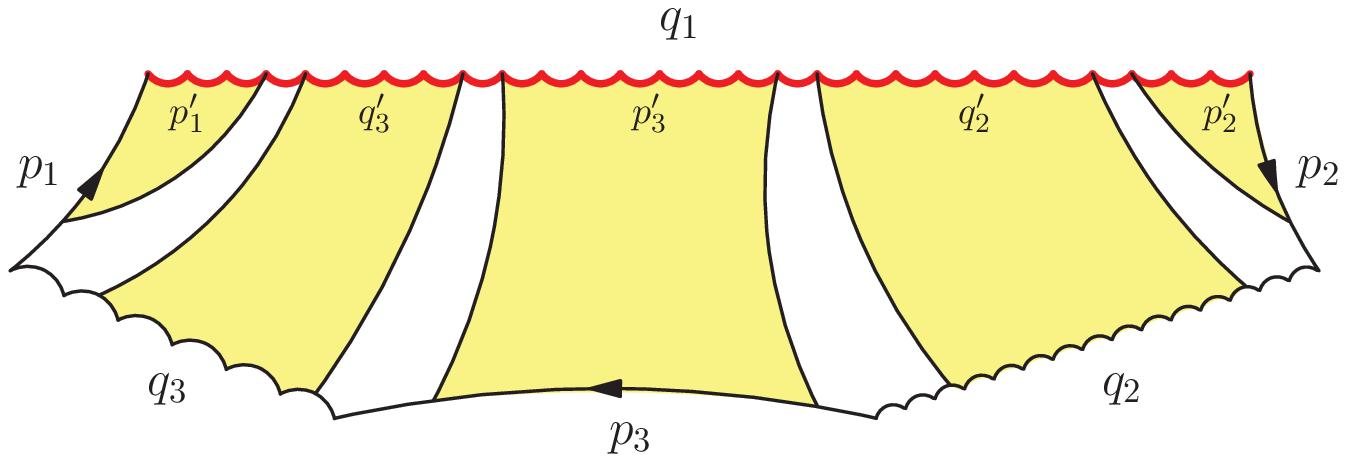}

\vspace*{-165mm}

\begin{center}
Fig. 5. The polygon $\mathcal{P}$.
\end{center}

\medskip

Then the path $q_1$ is covered by the union of its subpaths $p'_1,\dots,p'_n,q'_2,\dots,q'_n$
and at most $2n-2$ additional $g_1$-periods.
Therefore (and using the notation at the end of Section 4), we obtain
$$
N(q_1)\leqslant \overset{n}{\underset{i=1}{\sum}}\, N(p_i')+ \overset{n}{\underset{i=2}{\sum}}\, N(q_i')+2n-2.\eqno{(6.3)}
$$
We first estimate the numbers $N(p_i')$:
$$
\begin{array}{ll}
\displaystyle{N(p_i')= \frac{\ell(p_i')}{|g_1|_X}} & \displaystyle{\leqslant \frac{\varkappa\, d((p_i')_{-},(p_i')_{+})+\varepsilon}{|g_1|_X}}\vspace*{3mm}\\
& \displaystyle{\leqslant \frac{\varkappa (\ell(p_i)+2\nu)+\varepsilon}{|g_1|_X}}\vspace*{3mm}\\
& \displaystyle{=\frac{\varkappa (|a_i|_X+2\nu)+\varepsilon}{|g_1|_X}}.
\end{array}\eqno{(6.4)}
$$

\medskip




\noindent
In Claims 2 and 3 below we estimate the numbers $N(q_i')$.
Since the endpoints of $q_i'$ are
at distance at most $\nu$ from $q_i$, there exists a subpath $q_i''$ of $q_i$ or $\bar{q}_i$ such that
$d((q_i')_{-},(q_i'')_{-})\leqslant \nu$ and $d((q_i')_{+},(q_i'')_{+})\leqslant \nu$.
We need the following relation between $N(q_i')$ and $N(q_i'')$.

\medskip

{\bf Claim 1.} We have
$$
N(q_i'')\geqslant \frac{N(q_i')||g_1||_X}{|g_i|_X}-2\nu-2.\eqno{(6.5)}
$$

\medskip

{\it Proof.} The desired inequality follows from the following two estimations:
$$N(q_i'')\geqslant \frac{\ell(q_i'')}{|g_i|_X}-2,$$
$$
\ell(q_i'')\geqslant d((q_i'')_{-},(q_i'')_{+})\geqslant d((q_i')_{-},(q_i')_{+})-2\nu
=|g_1^{N(q_i')}|_X-2\nu
\geqslant N(q_i')||g_1||_X-2\nu.
$$
\hfill $\Box$

\medskip

In the following part of the proof we will use the function $f$ from Theorem~\ref{acylindric}. We set
$$
\alpha=\varkappa(2\nu+3+f(\nu)).\eqno{(6.6)}
$$

{\bf Claim 2.} If $g_1$ and $g_i$ are not commensurable, then
$$
N(q_i')\leqslant \alpha\frac{|g_i|_X}{|g_1|_X}+f(\nu).\eqno{(6.7)}
$$

\medskip

{\it Proof.} First consider the case $|g_1|_X\geqslant |g_i|_X$. Suppose that (6.7) is not valid. Then
$N(q_i')> f(\nu)$. Then, by Theorem~\ref{acylindric}, $g_1$ and $g_i$ are commensurable that contradicts
the assumption.

Now consider the case $|g_i|_X\geqslant |g_1|_X$. Suppose that (6.7) is not valid. Then
$$
N(q_i')>\alpha\frac{|g_i|_X}{|g_1|_X}.\eqno{(6.8)}
$$
Substituting (6.8) into (6.5), we deduce
$$
N(q_i'')\geqslant \alpha \frac{||g_1||_X}{|g_1|_X}-2\nu-2\,\overset{(2.1)}{\geqslant}\, \frac{\alpha}{\varkappa}-2\nu-2\,\overset{(6.7)}{>}\, f(\nu).
$$
By Theorem~\ref{acylindric} applied to $g_i$ and $g_1$, we obtain that these elements are commensurable.
A contradiction.\hfill $\Box$

\medskip
Now we set
$$
\beta=\varkappa(2\nu+3+F(\nu))L,\eqno{(6.9)}
$$
where $F$ is the function from Theorem~\ref{acylindric_1} and $L\geqslant 1$ is the constant from Lemma~\ref{elem_index}.

\medskip

{\bf Claim 3.} If $g_1$ and $g_i$ are commensurable, then
$$
N(q_i')\leqslant \beta\, {\bf Ind}_{[g_1]}([g_1]\vee [g_i]).\eqno{(6.10)}
$$

{\it Proof.} Suppose the converse, i.e.
$$
N(q_i')> \beta\, {\bf Ind}_{[g_1]}([g_1]\vee [g_i]).\eqno{(6.11)}
$$
Our nearest aim is to deduce the following two inequalities:
$$
N(q_i')> L\cdot {\bf Ind}_{[g_1]}([g_1]\vee [g_i])+F(\nu),\eqno{(6.12)}
$$
$$
N(q_i'')> L\cdot {\bf Ind}_{[g_i]}([g_i]\vee [g_1])+F(\nu).\eqno{(6.13)}
$$

The first inequality follows directly from the assumption (6.11) and the facts that $\beta \geqslant L+F(\nu)$ (since $\varkappa\geqslant 1$ in (6.9)) and ${\bf Ind}_{[g_1]}([g_1]\vee [g_i])\geqslant 1$. We prove the second one.
$$
\begin{array}{ll}
N(q_i'') & \displaystyle{\overset{(6.5)}{\geqslant} \frac{N(q_i')||g_1||_X}{|g_i|_X}}-2\nu-2\vspace*{3mm}\\
& \displaystyle{\overset{(6.11)}{\geqslant} \frac{\bigl(\beta \, {\bf Ind}_{[g_1]}([g_1]\vee [g_i])\bigr)||g_1||_X}{|g_i|_X}-2\nu-2}\vspace*{3mm}\\
& \displaystyle{\overset{(2.1)}{\geqslant} \frac{\bigl(\beta \, {\bf Ind}_{[g_1]}([g_1]\vee [g_i])\bigr)||g_1||_X}{\varkappa||g_i||_X}-2\nu-2}\vspace*{3mm}\\
& \displaystyle{\overset{(5.2)}{=} \frac{\beta}{\varkappa} {\bf Ind}_{[g_i]}([g_i]\vee [g_1])-2\nu-2}\vspace*{3mm}\\
& \overset{(6.9)}{\geqslant} L\cdot {\bf Ind}_{[g_i]}([g_i]\vee [g_1])+F(\nu).
\end{array}
$$
Thus, (6.12) and (6.13) are proved.
By Theorem~\ref{acylindric_1} and Lemma~\ref{nice_estimate}, there exist different phase vertices $x_1,x_2$ on $q_i'$ and different phase vertices $y_1,y_2$ on $q_i''$ such that $x_1^{-1}y_1=x_2^{-1}y_2$.
Then we can cut out a piece from $\mathcal{P}$ and glue the remaining pieces as shown in Figure~6.

\vspace*{-30mm}
\hspace*{-7mm}
\includegraphics[scale=0.75]{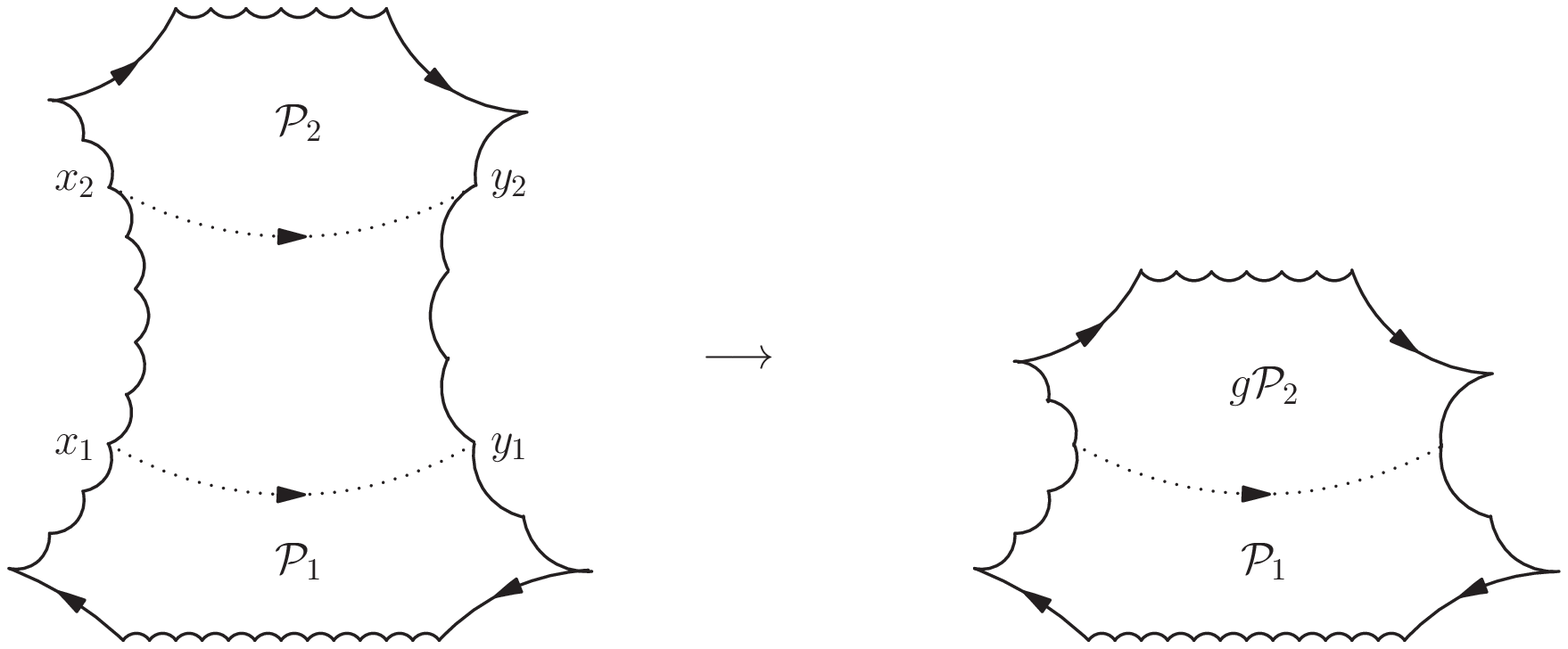}

\vspace*{-135mm}

\begin{center}
Fig. 6. Cutting out a piece from $\mathcal{P}$.
\end{center}

\medskip

More precisely, let $\mathcal{P}_1$ be the subpath of the (cyclic) path $\mathcal{P}$ from $y_1$ to $x_1$ and $\mathcal{P}_3$ be the subpath of $\mathcal{P}$ from $x_2$ to $y_2$. We consider the polygon $\mathcal{P}'$
obtained by gluing the endpoints of $\mathcal{P}_1$ to the corresponding endpoints of the left translation $g\mathcal{P}_3$, where $g=x_1x_2^{-1}=y_1y_2^{-1}$.
The new polygon $\mathcal{P}'$ corresponds to a solution of (6.1) with smaller value $|k_1|+\dots +|k_n|$.
A contradiction.
\hfill $\Box$

\medskip

Thus, the summands in (6.3) are estimated in (6.4) and in Claims 2 and 3.
This proves the inequality (6.2) for some universal constant $M$.
\hfill $\Box$

\section{Theorem A and its proof}

\noindent
{\bf Theorem A.} {\it Let $G$ be an acylindrically hyperbolic group with respect to a generating set $X$.
Then there exists a constant $M>1$
such that for any exponential equation
$$
a_1g_1^{x_1}a_2g_2^{x_2}\dots a_ng_n^{x_n}=1\eqno{(7.1)}
$$
with constants $a_1,g_1,\dots,a_n,g_n$ from $G$ and variables $x_1,\dots,x_n$, if this equation is solvable over $\mathbb{Z}$, then there exists a solution $(k_1,\dots,k_n)$ with
$$
|k_j|\leqslant \Bigl(n^2+\overset{n}{\underset{i=1}\sum}\,\frac{|a_i|_X}{|g_j'|_X}+
\underset{g_i\notin {\text{\rm Com}}(g_j)}\sum\,\frac{|g_i|_X}{|g_j'|_X}+
\underset{g_i\in {\text{\rm Com}}(g_j)}{\sum}\, |g_i|_X
\Bigr)\cdot M\eqno{(7.2)}
$$
for all $j$ corresponding to loxodromic $g_j$; here $g'_j$ is an element shortest
in the conjugacy class of $g_j$ with respect to $X$.

This implies that if the equation (7.1) is solvable over $\mathbb{Z}$, then there exists a solution $(k_1,\dots,k_n)$ with the universal estimation
$$
|k_j|\leqslant \Bigl(n^2+\overset{n}{\underset{i=1}\sum}\,|a_i|_X+\overset{n}{\underset{i=1}{\sum}}\, |g_i|_X\Bigr)\cdot M\eqno{(7.3)}
$$
for all $j$ corresponding to loxodromic $g_j$.
}

\medskip

{\it Proof.}
We first consider the special case, where all $g_i$ are loxodromic.
By Lemma~\ref{loxo}, for each $g_i$, there exists $h_i\in G$ such that
$g_i'=h_i^{-1}g_ih_i$ has minimal length in the conjugacy class of $g_i$ and
$$
|h_i|_X\leqslant K|g_i|_X.\eqno{(7.4)}
$$
For convenience we set $g_0=g_n$ and $h_0=h_n$. Now we rewrite (7.1) as
$$
a_1'(g_1')^{x_1}a_2'(g_2')^{x_2}\dots a_n'(g_n')^{x_n}=1,
$$
where $a_i'=h_{i-1}^{-1}a_ih_i$ for $i=1,\dots,n$. Since $g_i'$ are loxodromic and shortest in their conjugacy classes, by Lemma~\ref{first_main lemma}, there exists a solution $(k_1,\dots,k_n)$ of equation (7.1) with
$$
|k_j|\leqslant \Bigl(n^2+\overset{n}{\underset{i=1}\sum}\,\frac{|a'_i|_X}{|g'_j|_X}+
\underset{g_i\notin {\text{\rm Com}}(g_j)}\sum\,\frac{|g'_i|_X}{|g'_j|_X}+
\underset{g_i\in {\text{\rm Com}}(g_j)}{\sum}
{\bf Ind}_{[g_j]}([g_j]\vee [g_i])
\Bigr)\cdot M_1\eqno{(7.5)}
$$
for all $j=1,\dots,n$, where $M_1$ is a universal constant.
We set
$$
M_2=2M_1K \frac{L^2}{{\bf inj}(G,X)}.
$$
Then (7.2) with $M=M_2$ follows from (7.5) with the help of the following claim.

{\bf Claim.} We have

1) $
|a_i'|_X\leqslant |a_i|_X+K\bigl(|g_{i-1}|_X+|g_i|_X\bigr).
$

2)
$
|g_i'|_X\leqslant |g_i|_X.
$

3) $$
{\bf Ind}_{[g_j]}([g_j]\vee [g_i])
\leqslant  \frac{L^2}{{\bf inj}(G,X)}\cdot |g_i|_X.
$$
{\it Proof.} The first statement follows from the definition of $a_i'$ and (7.4),
the second from the definition of $g_i'$, and the third from Lemma~\ref{estim_k}.\hfill $\Box$

\medskip

Now we consider the general case.
Let $\mathcal{E}$ (resp. $\mathcal{L}$) be the set of the indexes $i\in \{1,\dots,n\}$ for which $g_i$ is elliptic (resp. loxodromic). We have $\mathcal{E}\cup \mathcal{L}=\{1,\dots,n\}$.
Note that by Lemma~\ref{elli}, if $i\in \mathcal{E}$, then
$$
\begin{array}{ll}
|g_i^{x_i}|_X & \leqslant 2K|g_i|_X +(8\delta+1)\vspace*{2mm}\\
 & \leqslant 2K(8\delta+2)|g_i|_X.
\end{array}\eqno{(7.6)}
$$
for any choice of $x_i\in \mathbb{Z}$.
For any two consecutive numbers $s,t\in \mathcal{L}$, let $b_t$ be the product of the factors in (7.1)
between $g_s^{x_s}$ and $g_t^{x_t}$, i.e.
$$
b_t=a_{s+1}g_{s+1}^{x_{s+1}}\dots g_{t-1}^{x_{t-1}}a_t.\eqno{(7.7)}
$$
Then we can reduce to the considered  case (all $g_i$ are loxodromic) by writing
$$
a_1g_1^{x_1}\dots a_ng_n^{x_n}=\underset{i\in \mathcal{L}}{\prod}b_ig_i^{x_i}.
$$
From this case we have
$$
|k_j|\leqslant \Bigl(n^2+\frac{1}{|g_j'|_X}\Bigl(\overset{n}{\underset{i\in \mathcal{L}}\sum}\,|b_i|_X+
\underset{i\in \mathcal{L}, g_i\notin {\text{\rm Com}}(g_j)}\sum\,|g_i|_X\Bigr) +
\underset{i\in \mathcal{L}, g_i\in {\text{\rm Com}}(g_j)}{\sum}|g_i|_X
\Bigr)\cdot M_2\eqno{(7.8)}
$$
for any $j\in \mathcal{L}$ and some universal constant $M_2$.
Now we estimate the sums in the internal brackets.
First observe that for any choice of $x_i$, we have
$$
\begin{array}{ll}
\underset{i\in \mathcal{L}}{\sum}|b_i|_X & \overset{(7.7)}{\leqslant}\, \overset{n}{\underset{i=1}{\sum}}|a_i|_X+\underset{i\in \mathcal{E}}{\sum}|g_i^{x_i}|_X.\vspace*{2mm}\\
& \overset{(7.6)}{\leqslant}\, \overset{n}{\underset{i=1}{\sum}}|a_i|_X+2K(8\delta+2)\underset{i\in \mathcal{E}}{\sum}|g_i|_X\vspace*{2mm}\\
& \leqslant\, \overset{n}{\underset{i=1}{\sum}}|a_i|_X+2K(8\delta+2)\underset{i\in \mathcal{E}, g_i\notin {\text{\rm Com}}(g_j)}{\sum}|g_i|_X\
\end{array}
$$
The last inequality is satisfied since the condition $g_i\notin {\text{\rm Com}}(g_j)$ is automatically satisfied for $i\in \mathcal{E}$ (elliptic and loxodromic elements are not commensurable). Therefore the sum in the internal brackets in (7.8) does not exceed
$$
2K(8\delta+2)\Bigl(\overset{n}{\underset{i=1}{\sum}}\,|a_i|_X+
\underset{g_i\notin {\text{\rm Com}}(g_j)}\sum\,|g_i|_X\Bigr).
$$
Then (7.2) is satisfied for $M=2M_2K(8\delta+2)$.
The last statement of the main theorem follows from (7.2) and $|g_j'|_X\geqslant 1$.
\hfill $\Box$



\section{Theorems B and C and their proofs}

In subsection 8.1 we recall some definitions and statements about hyperbolically embedded subgroups and weakly hyperbolic groups.
Theorem B is formulated and proved in subsection 8.2.
Theorem C is deduced from Theorems A$'$ and B in subsection 8.3.

\subsection{Some definitions and statements from~\cite{DOG}}
Let $G$ be a group, $\{H_{\lambda}\}_{\lambda\in \Lambda}$ a collection of subgroups of $G$.
A subset $X$ of $G$ is called a {\it relative generating set of $G$ with respect to}
$\{H_{\lambda}\}_{\lambda\in \Lambda}$ if $G$ is generated by $X$ together with the union of all $H_{\lambda}$.
All relative generating sets are assumed to be symmetric.
We define
$$
\mathcal{H}=\bigsqcup_{\lambda\in\Lambda}H_{\lambda}.
$$


In this section, we always assume that $X$ is a relative generating set of $G$ with respect to $\{H_{\lambda}\}_{\lambda\in \Lambda}$.

\begin{definition} {\rm (see~\cite[Definition 4.1]{DOG})
The group $G$ is called {\it weakly hyperbolic} relative to $X$ and $\{H_{\lambda}\}_{\lambda\in \Lambda}$ if the Cayley graph $\Gamma(G, X\sqcup \mathcal{H})$ is hyperbolic.
}
\end{definition}

\noindent
We consider the Cayley graph
$\Gamma(H_{\lambda},H_{\lambda})$ as a complete subgraph of $\Gamma(G,X\sqcup \mathcal{H})$.

\begin{definition}\label{definition_relative_metrics}
{\rm (see~\cite[Definition 4.2]{DOG})
For every $\lambda\in \Lambda$, we introduce a {\it relative metric}
$\widehat{d}_{\lambda}:H_{\lambda}\times H_{\lambda}\rightarrow [0,+\infty]$ as follows:

Let $a,b\in H_{\lambda}$. A path
in $\Gamma(G,X\sqcup \mathcal{H})$ from $a$ to $b$ is called {\it $H_{\lambda}$-admissible} if it has no edges in the subgraph $\Gamma(H_{\lambda},H_{\lambda})$.

The distance $\widehat{d}_{\lambda}(a,b)$ is defined to be the length of a shortest
{\it $H_{\lambda}$-admissible} path connecting $a$ to $b$ if such exists.
If no such path exists, we set $\widehat{d}_{\lambda}(a,b)=\!
\infty$.


}
\end{definition}

\begin{definition}\label{def_hyperb_embedd} {\rm (see~\cite[Definition 4.25]{DOG})
Let $G$ be a group, $X$ a symmetric subset of $G$. A collection of subgroups $\{H_{\lambda}\}_{\lambda\in \Lambda}$
of $G$ is called {\it hyperbolically embedded in $G$ with respect to $X$}
(we write $\{H_{\lambda}\}_{\lambda\in \Lambda}\hookrightarrow_h (G,X)$) if the following hold.

\begin{enumerate}
\item[(a)] The group $G$ is generated by $X$ together with the union of all $H_{\lambda}$ and the Cayley graph
$\Gamma(G,X\sqcup \mathcal{H})$ is hyperbolic.

\item[(b)] For every $\lambda\in \Lambda$, the metric space $(H_{\lambda},\widehat{d}_{\lambda})$ is
proper. That is, any ball of finite radius in $H_{\lambda}$ contains finitely many elements.
\end{enumerate}


}
\end{definition}




\begin{definition}\label{components}
{\rm (see~\cite[Definition 4.5]{DOG})
Let $q$ be a path in the Cayley graph $\Gamma(G,X\sqcup \mathcal{H})$. A non-trivial subpath $p$ of $q$
is called an {\it $H_{\lambda}$-subpath}, if the label of $p$ is a word in the alphabet $H_{\lambda}$.
An $H_{\lambda}$-subpath $p$ of $q$ is called an {\it $H_{\lambda}$-component} if $p$ is not contained in a longer subpath of $q$ with this property. Two $H_{\lambda}$-components $p_1,p_2$ of a path $q$ in $\Gamma(G,X\sqcup \mathcal{H})$ are called {\it connected} if there exists a path $\gamma$ in $\Gamma(G,X\sqcup \mathcal{H})$ that connects some vertex of $p_1$ to some vertex of $p_2$, and ${\bf Lab}(\gamma)$ is a word consisting only of letters from
$H_{\lambda}$.

Note that we can always assume that $\gamma$ has length at most 1 as every element of $H_{\lambda}$
is included in the set of generators. An $H_{\lambda}$-component $p$ of a path $q$ in $\Gamma(G,X\sqcup \mathcal{H})$ is {\it isolated} if it is not connected to any other component of $q$.

}
\end{definition}

Given a path $p$ in $\Gamma(G, X\sqcup \mathcal{H})$, the canonical image of ${\bold{Lab}}(p)$ in $G$ is denoted by ${\bold {Lab}}_G(p)$.



\begin{definition}\label{n-gon} {\rm (see~\cite[Definition 4.13]{DOG})
Let $\varkappa \geqslant 1$, $\varepsilon\geqslant 0$, and $m\geqslant 2$. Let $\mathcal{P}=p_1\dots p_m$ be an $m$-gon in $\Gamma(G, X\sqcup \mathcal{H})$ and let $I$ be a subset of the set of its sides $\{p_1,\dots,p_m\}$
such that:

1) Each side $p_i\in I$ is an isolated $H_{\lambda_i}$-component of $\mathcal{P}$ for some $\lambda_i\in \Lambda$.

2) Each side $p_i\notin I$ is a $(\varkappa,\varepsilon)$-quasi-geodesic.

\medskip
\noindent
We denote $s(\mathcal{P},I)=\underset{p_i\in I}{\sum} \widehat{d}_{\lambda_i}(1,{\bold{Lab}}_G(p_i))$.

}
\end{definition}

\begin{proposition}\label{Proposition_Osin}
{\rm (see~\cite[Proposition 4.14]{DOG})}
Suppose that $G$ is weakly hyperbolic relative to $X$ and $\{H_{\lambda}\}_{\lambda\in \Lambda}$.
Then for any $\varkappa\geqslant 1$, $\varepsilon\geqslant 0$, there exists a constant $C(\varkappa,\varepsilon)>0$ such that
for any $m$-gon $\mathcal{P}$ in $\Gamma(G,X\sqcup \mathcal{H})$ and any subset $I$ of the set of its sides satisfying conditions
of Definition~\ref{n-gon}, we have $s(\mathcal{P},I)\leqslant C(\varkappa,\varepsilon)m$.
\end{proposition}

\subsection{Elliptical exponential equations over a group with given hyperbolically embedded subgroups}\hspace*{-3mm}.

\medskip
\noindent
{\bf Theorem B.} {\it Let $G$ be a group, $\{H_{\lambda}\}_{\lambda\in \Lambda}$ a collection of subgroups of $G$, and $X$ a symmetric relative
generating set of $G$ with respect to $\{H_{\lambda}\}_{\lambda\in \Lambda}$.
Suppose that $\{H_{\lambda}\}_{\lambda\in \Lambda}$ is hyperbolically embedded in $G$ with respect to $X$.
Then any exponential equation
$$
a_1g_1^{x_1}a_2g_2^{x_2}\dots a_ng_n^{x_n}=1\eqno{(8.1)}
$$
with $a_1,\dots,a_n\in G$  and $g_1,\dots,g_n\in \mathcal{H}=\bigsqcup_{\lambda\in\Lambda}H_{\lambda}$
is equivalent to a finite disjunction of finite systems of equations,
$$
\overset{k}{\underset{i=1}{\bigvee}}\overset{\ell_i}{\underset{j=1}{\bigwedge}} E_{ij},
$$
such that
\begin{enumerate}
\item[{\rm (1)}] each $E_{ij}$ is  an exponential equation over some $H_{\lambda}$,
or a trivial equation of kind $g_{ij}=1$, where $g_{ij}$ is an element of $G$,

\item[{\rm (2)}] for any $i=1,\dots, k$, the sets of variables of
$E_{i,j_1}$ and $E_{i,j_2}$ are disjoint if $j_1\neq j_2$.

\end{enumerate}

Let $\Lambda_0=\{\lambda_1,\dots,\lambda_n\}$ be a subset of $\Lambda$ such that $g_i\in H_{\lambda_i}$, $i=1,\dots,n$, and let
$L=n+\overset{n}{\underset{i=1}{\sum}} |a_i|_{X\cup \mathcal{H}}$.
Then these systems of equations can be algorithmically written if for any $\lambda\in \Lambda_0$,
there is an algorithm computing the following finite subsets of~$H_{\lambda}$:
$$
H_{\lambda,L}=\{h\in H_{\lambda}\,|\, \widehat{d}_{\lambda}(1,h)\leqslant C(1,1)\cdot L\},\eqno{(8.2)}
$$
where $C(1,1)$
is the constant from Proposition~\ref{Proposition_Osin}.}



\bigskip

{\it Proof.}
To describe the desired family of systems of equations formally, we first introduce definitions (a)-(b) below.
Let $f:\{1,\dots,n\}\rightarrow \Lambda_{0}$ be a map such that $g_i\in H_{f(i)}$ for $i=1,\dots,n$.
For any $\lambda\in \Lambda_0$ we define the set
$$
H_{\lambda}^{\ast}=H_{\lambda}\cup \{g_i^{x_i}\,|\, f(i)=\lambda\},
$$
where $g_i^{x_i}$ is considered as a single letter.
We also define $\mathcal{H}^{\ast}=\underset{\lambda\in \Lambda_0}{\bigsqcup} H_{\lambda}^{\ast}$. Thus,
$\mathcal{H}^{\ast}=\mathcal{H}\cup \{g_1^{x_1},\dots,g_n^{x_n}\}$.

We represent each element $a_i$ by a word $A_i$ (not necessarily of minimal possible length) in the alphabet
$X\sqcup \mathcal{H}$.
Then the expression on the left side of (8.1) can be represented by the word $\bold{W}=A_1g_1^{x_1}A_2g_2^{x_2}\dots A_ng_n^{x_n}$ in the alphabet $X\sqcup \mathcal{H}^{\ast}$. Let $L$ be the length of this word; we have $L=n+\overset{n}{\underset{i=1}{\sum}}\, |A_i|$.

We consider a closed disc $D$ such that its oriented boundary $\partial D$ is divided into $L$ consecutive paths
$s_1,s_2,\dots,s_L$ labelled by the elements of $X\sqcup \mathcal{H}^{\ast}$
so that the label of $\partial D$ coincides with the cyclic word $\bold{W}$.
Thus the equation (8.1) can be written in the form ${\text{\bf Lab}}(\partial D)=1$.


\medskip

(a) Let $\lambda\in \Lambda_0$ and let $P$ be a nontrivial subpath of the cyclic combinatorial path $\partial D=s_1s_2\dots s_L$. The subpath $P$ is called an {\it $H_{\lambda}^{\ast}$-subpath} of $\partial D$ if the label of $P$ is a word in the alphabet $H_{\lambda}^{\ast}$.
An $H_{\lambda}^{\ast}$-subpath $P$ of $\partial D$ is called an {\it $H_{\lambda}^{\ast}$-component} if $P$ is not contained in a longer subpath of $\partial{D}$ with this property.
Sometimes we will skip the subscript $\lambda$ and call $P$ an {\it $\mathcal{H}^{\ast}$-component} of $\partial D$.

The cyclic combinatorial path $\partial D$ can be written as $\partial D=Q_1P_1\dots Q_rP_r$, where $P_1,\dots,P_r$ are all $\mathcal{H}^{\ast}$-components of $\partial D$.
We say that an $\mathcal{H}^{\ast}$-component $P$ is {\it special}
if the label of $P$ contains the letter $g_j^{x_j}$ for some $j\in \{1,\dots,n\}$.

\medskip


(b) A region $R$ in $D$ homeomorphic to a closed disc is called an {\it $H_{\lambda}^{\ast}$-region} if its boundary has the form $U_1E_1\dots U_sE_s$, where
$U_1,\dots ,U_s$ are $H_{\lambda}^{\ast}$-components for the same $\lambda\in \Lambda_0$ and at least one of them is special,
and $E_1,\dots, E_s$ are simple paths in $D$ whose interiors lie in the interior of $D$, see Figure 7.
We call these paths {\it internal sides} of $R$.
We say that the {\it internal sides of $R$ are boundedly  labelled}, if each $E_i$ is labelled by an element of
$H_{\lambda,L}$, see (8.2).

\vspace*{-20mm}
\hspace*{10mm}
\includegraphics[scale=0.5]{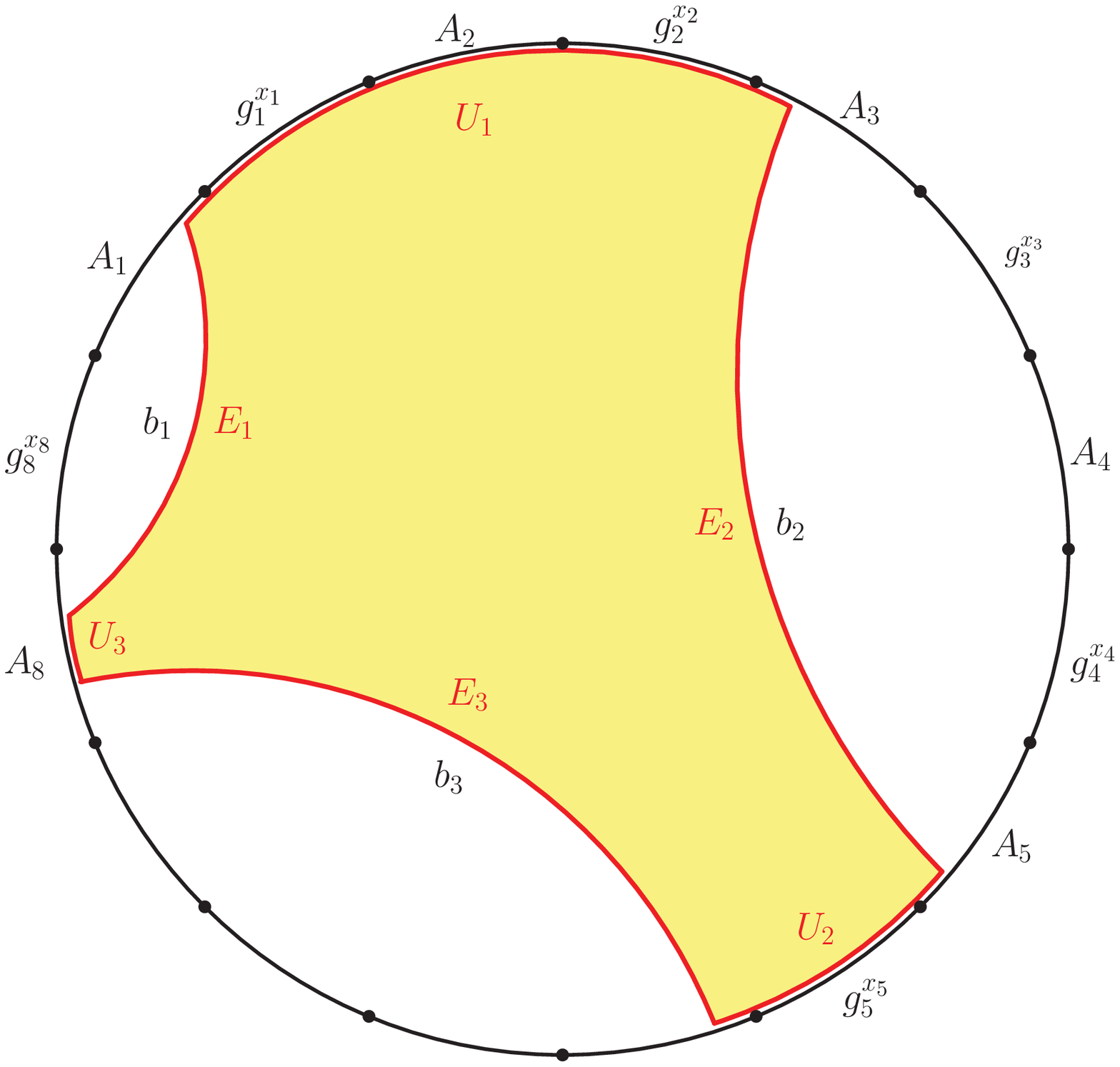}

\vspace*{-40mm}

\begin{center}
Fig.7. An example of an $H_{\lambda}^{\ast}$-region, where $U_1$ and $U_2$ are special $\mathcal{H}^{\ast}$-components.
\end{center}

Note that the set $H_{\lambda, L}$ is finite, since the metric space $(H_{\lambda}, \widehat{d}_{\lambda})$ is locally finite.

A collection of regions $\mathcal{R}=\{R_1,\dots ,R_t\}$, where each $R_i$ is an $H_{\lambda (i)}^{\ast}$-region for some $\lambda(i)\in \Lambda_0$ is called {\it admissible} if
the intersection of $R_i$ and $R_j$ is either empty or consists of one or two points on $\partial D$.
We do not distinguish two admissible collections of regions $\mathcal{R}$ and $\mathcal{R}'$ if there exists an isotopy of $D$ fixing $\partial D$ and carrying the elements of $\mathcal{R}$ to the elements of $\mathcal{R}'$.

A collection of regions $\mathcal{R}=\{R_1,\dots ,R_t\}$ is called {\it complete} if it is admissible and any
special $\mathcal{H}^{\ast}$-component $P_i$ is contained in the boundary of some $R_j\in \mathcal{R}$.

\medskip




\medskip

(c) Let $\mathcal{R}=\{R_1,\dots,R_s\}$ be any complete collection of regions with boundedly labelled internal sides.
Let $R_{s+1},\dots,R_{s+t}$ be the components of the closure of $D\setminus \cup \,\mathcal{R}$.
Then $\mathcal{R}$ determines a system of exponential equations over $G$, namely
$$
{\text{\bf Eq}}(\mathcal{R}):
\begin{cases}
{\text{\bf Lab}}(\partial R_1)=1, &\vspace*{2mm}\\
\dots &\vspace*{2mm}\\
{\text{\bf Lab}}(\partial R_{s+t})=1.
\end{cases}
$$

Clearly, any solution of the system ${\text{\bf Eq}}(\mathcal{R})$ satisfies the equation ${\text{\bf Lab}}(\partial{D})=1$. Moreover, the first $s$ equations of this system are exponential equations over $H_{\lambda}$, where $\lambda$ goes through $\Lambda_0$,  The last $t$ equations have the form $U=1$, where $U$ is a word in the alphabet $X\sqcup \mathcal{H}$  (i.e. it has no letters $g_i^{x_i}$).

\medskip

{\bf Claim.} The set $\mathfrak{F}$ of all complete collections of regions with boundedly labelled internal sides is finite. Each solution of (8.1) satisfies the system ${\text{\bf Eq}}(\mathcal{R})$ for some $\mathcal{R}\in \mathfrak{F}$.

\medskip

{\it Proof.}
The finiteness of $\mathfrak{F}$ follows from the finiteness of $H_{\lambda,L}$ for any $\lambda\in \Lambda_0$.

Suppose that $\overline{k}=(k_1,\dots,k_n)$ is some solution of the equation (8.1). For brevity,
we introduce the following two definitions.

\medskip

{\it Definition 1.} 
Let $\Delta$ be a graph with  edges labelled by elements of the alphabet $X\sqcup \mathcal{H}^{\ast}$.
A graph map $\psi: \Delta \rightarrow \Gamma(G,X\sqcup \mathcal{H})$ is called a {\it $\overline{k}$-map}, if $\psi$ maps edges labelled by elements of $X\sqcup \mathcal{H}$ to edges labelled by the same elements,
and edges labelled by $g_i^{x_i}$ to edges labelled by~$g_i^{k_i}$.

\medskip

{\it Definition 2.} Let $\mathcal{R}$ be an admissible collection of regions in $D$. We denote by $D_{\mathcal{R}}$ the $CW$-complex obtained from $D$ by subdivision along all internal sides of all regions from $\mathcal{R}$.  We use the notation $D_{\mathcal{R}}^{(1)}$ for the graph associated with the 1-skeleton
of $D_{\mathcal{R}}$. Thus, the edges of $D_{\mathcal{R}}^{(1)}$ are the paths $s_1,\dots,s_L$ and the internal sides of all regions from $\mathcal{R}$.

\medskip

Observe that the above claim can be directly deduced from the following statement.

\medskip

{\it Statement.} Let $\overline{k}=(k_1,\dots,k_n)$ be an arbitrary solution of the equation (8.1) and let
$\mathcal{P}$ be some closed path in $\Gamma(G,X\sqcup \mathcal{H})$
with the label $A_1g_1^{k_1}\dots A_ng_n^{k_n}$.
 Then there exists a complete collection $\mathcal{R}$ of regions in $D$ with boundedly labelled
internal sides such that the $\overline{k}$-map $\partial D\rightarrow \mathcal{P}$ extends
to a $\overline{k}$-map $D_{\mathcal{R}}^{(1)} \rightarrow \Gamma(G,X\sqcup \mathcal{H})$.

\medskip

It remains to prove this statement.
Note that if $p$ is an arbitrary $H_{\lambda}$-component of $\mathcal{P}$, then there exists an edge $e$ in $\Gamma(G,X\cup \mathcal{H})$ such that $pe$ is a closed path in  $\Gamma(G,X\sqcup \mathcal{H})$,
and we have ${\text{\bf Lab}}(e)\in H_{\lambda}$.

Let $p_{j_1}e_1p_{j_2}e_2\dots p_{j_m}e_m$ be a closed path in $\Gamma(G,X\sqcup \mathcal{H})$ such that
$p_{j_1},p_{j_2},\dots p_{j_m}$ are  $H_{\lambda}$-components of $\mathcal{P}$ for the same $\lambda\in \Lambda_0$,
$j_1<j_2<\dots <j_m$ (where we use the cyclic ordering on $\mathbb{Z}_r$), the corresponding component $P_{j_1}$ of $\partial D$ is special, $e_1,e_2,\dots, e_m$ are edges in $\Gamma(G,X\sqcup \mathcal{H})$ with labels from $H_{\lambda}$, and $m$ is maximal with these properties.

Let $q_i$ be a subpath in $\mathcal{P}$ such that $(q_i)_{-}=(e_i)_{-}$, $(q_i)_{+}=(e_i)_{+}$, $i=1,\dots, m$. Denote $\mathcal{P}_i=q_ie_i^{-1}$.
We claim that $e_i^{-1}$ is an isolated $H_{\lambda}$-component in $\mathcal{P}_i$ for any~$i$.
Indeed, since $p_{j_i}$ and $p_{j_{i+1}}$ are $H_{\lambda}$-components of $\mathcal{P}$,
the labels of the first and the last edges of $q_i$ do not lie in $H_{\lambda}$. Therefore
$e_i^{-1}$ is an $H_{\lambda}$-component in $\mathcal{P}_i$. Since $m$ is maximal, this component is isolated in $\mathcal{P}_i$. By Proposition~\ref{Proposition_Osin}, we have ${\text{\bf Lab}}(e_i)=b_i$ for some  $b_i\in H_{\lambda,L}$.

Now we lift the edges $e_i$ to $D$, i.e., for any $e_i$ let $E_i$ be the directed chord in $D$ such that
$(E_i)_{-}=(P_{j_i})_{+}$, $(E_i)_{+}=(P_{j_{i+1}})_{-}$; we set ${\text{\bf Lab}}(E_i)=b_i$, see Figure~8.
Let $R$ be the $H_{\lambda}^{\ast}$-region in $D$ with the boundary $\partial R=P_{j_1}E_1\dots P_{j_m}E_m$.
Let $D_i$ be the closure of the component of $D\setminus R$, which contains $E_i$ in its boundary, $i=1,\dots,m$.
By induction, there exists a complete collection $\mathcal{R}_i$ of regions in $D_i$
with boundedly labelled internal sides such that the $\overline{k}$-map
$\partial D_i\rightarrow \mathcal{P}_i$ extends to a $\overline{k}$-map
$(D_i)_{\mathcal{R}_i}^{(1)}\rightarrow \Gamma(G,X\sqcup \mathcal{H})$.
Then the collection $\mathcal{R}=\{R\}\bigcup \overset{m}{\underset{i=1}{\cup}}\,\mathcal{R}_i$
satisfies the above statement.
\hfill $\Box$ $\Box$

\vspace*{-15mm}
\hspace*{10mm}
\includegraphics[scale=0.5]{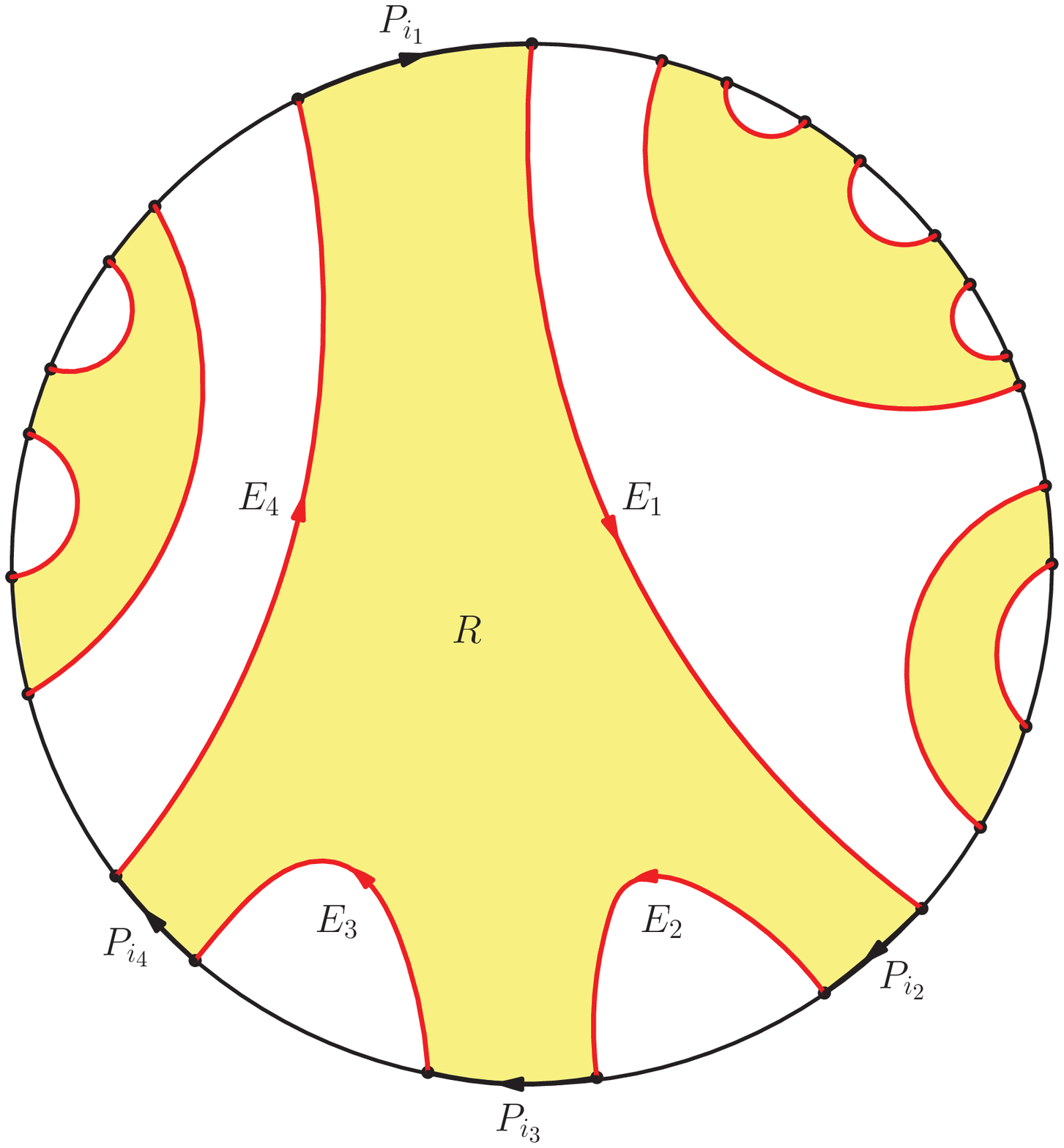}

\vspace*{-35mm}

\begin{center}
Fig. 8. Illustration to the proof of the statement.
\end{center}

\bigskip

\subsection{Proof of Theorem C}
We first prove two auxiliary lemmas about relatively hyperbolic groups, which have algorithmic character.
We rely on the manuscript of Osin~\cite{Osin_0}.

\begin{remark} Let $G$ be a group relatively hyperbolic with respect to a collection of subgroups $\{H_{\lambda}\}_{\lambda\in \Lambda}$, and let $X$ be a finite relative generating set of $G$.
It is well known that any element of $G$ has exactly one of the following three types:
(1) parabolic, (2) non-parabolic of finite order, (3) loxodromic with respect to $X\cup \mathcal{H}$.
\end{remark}

\begin{lemma}\label{efficient_type}
Let $G$ be a group which is relatively hyperbolic with respect to a finite collection of its subgroups $\mathbb{H}=\{H_1,\dots,H_m\}$. Suppose that
\begin{enumerate}
\item[\rm{(a)}] $G$ is finitely generated,

\item[\rm{(b)}] each subgroup $H_i$ is given by a recursive presentation and has solvable word problem,

\item[\rm{(c)}] $G$ is given by a finite relative presentation $\mathcal{P}=\langle X\,|\, \mathcal{R}\rangle$  with respect to  $\mathbb{H}$, where $X$ is a finite set generating $G$,

\item[\rm{(d)}] the hyperbolicity constant $\delta$ of the Cayley graph $\Gamma(G,X\cup \mathcal{H})$ is known.
\end{enumerate}

Then the question about the type of an element $g\in G$ (given as a word in the alphabet $X\sqcup \mathcal{H}$)
is algorithmically decidable.


\end{lemma}

{\it Proof.} By~\cite[Theorem 5.6]{Osin_0}, we determine whether $g$ is parabolic or not.
Suppose that $g$ is nonparabolic. We show how to determine whether the order of $g$ is finite or not.

By~Lemma 4.5 from~\cite{Osin_0} (together with the last line of its proof) combined with Corollary~4.4 rom~\cite{Osin_0}, any element of finite order in $G$ is conjugate to an element of the set
$$
S=\{a\in G\,|\, |a|_X\leqslant B\cdot (8\delta+1)^2\},
$$
where $B=2C\, \underset{R\in \mathcal{R}}{\max} |R|_{X\cup \mathcal{H}}$, and $C$ is the constant
in the relative Dehn function $D_G^{rel}$.
Since $X$ is finite, we can find the set $S$ efficiently.
Let $I=\{0,1,\dots,|S|\}$,
For $i\in I$ we check whether $g^i$ is conjugate to an element of $S$, see Theorem 5.13 from~\cite{Osin_0}. If for some $i\in I$ the element $g^i$ is not conjugate to an element of $S$, then $g^i$ (and hence $g$) is loxodromic.
If every element $g^i$, $i\in I$, is conjugate to an element of $S$, then there exist two different numbers $i,j\in I$ such that $g^i$ is conjugate to $g^j$. In this case $g$ cannot be loxodromic, hence $g$ has a finite order.\hfill $\Box$

\medskip

\begin{lemma}\label{efficient_M}
Let $G$ be a finitely generated group which is relatively hyperbolic with respect to a finite collection of its subgroups $\{H_1,\dots,H_m\}$.
Suppose that $G$ is given by a finite relative presentation $\mathcal{P}=\langle X\,|\, \mathcal{R}\rangle$  with respect to  $\{H_1,\dots,H_m\}$, where $X$ is a finite set generating $G$. Suppose we know the hyperbolicity constant $\delta$ of the Cayley graph $\Gamma(G,X\sqcup \mathcal{H})$.
Then the constant $M$ from Theorem $A$ can be algorithmically computed.
\end{lemma}


{\it Proof.} We may assume that all subgroups $H_i$ are proper.
Then, by Proposition 5.2 from~\cite{Osin_1}, $G$ is acylindrically hyperbolic with respect to $X\cup \mathcal{H}$.
We claim that the following functions and constants can be computed in terms of $|X|$, $\delta$, and $\underset{r\in \mathcal{R}}{\max}|r|_{X\cup \mathcal{H}}$:

\medskip

$\bullet$ the functions $R$ and $N$ from Definition~\ref{acyl_action},

$\bullet$ the constant $L$ from Lemma~\ref{elem_index},

$\bullet$ the injectivity radius ${\bf inj}(G,X\cup \mathcal{H})$, see the paragraph before Definition~\ref{L}.

\medskip

Indeed, by the proof of Proposition 5.2 from~\cite{Osin_1}, one can take $R(\varepsilon)= 6\varepsilon +2$, $N(\varepsilon)=(6\varepsilon+2)|B_X(2\varepsilon)|$. By the proof of Lemma~6.8 from~\cite{Osin_1}, one can compute $L$ in terms of $\delta$ with the help of the functions $R$ and $N$. Finally, one can compute
${\bf inj}(G,X\cup \mathcal{H})$ in terms of $|X|$, $\delta$ and $\underset{r\in \mathcal{R}}{\max}|r|_{X\cup \mathcal{H}}$,
see the proof of Theorem 4.25 from~\cite{Osin_0}.

Following the proof of Theorem~A$'$, where these functions and constants were used, one can compute $M$. \hfill $\Box$

\noindent
{\bf Theorem C.}
{\it Let $G$ be a group relatively hyperbolic with respect to a finite collection of
subgroups $\{H_1,\dots,H_m\}$.
Suppose that $G$ is finitely generated, each subgroup $H_i$ is given by a recursive presentation and has solvable word problem,
$G$ is given by a finite relative presentation $\mathcal{P}=\langle X\,|\, \mathcal{R}\rangle$  with respect to  $\{H_1,\dots,H_m\}$, where $X$ is a finite set generating $G$, and that the hyperbolicity constant $\delta$ of the Cayley graph $\Gamma(G,X\cup \mathcal{H})$ is known, $\mathcal{H}=\overset{m}{\underset{i=1}{\bigsqcup}}H_i$.

Then there exists an algorithm which for any exponential equation $E$ over $G$ finds a finite disjunction $\Phi$ of finite systems of equations,
$$
\Phi:=\overset{k}{\underset{i=1}{\bigvee}}\overset{\ell_i}{\underset{j=1}{\bigwedge}} E_{ij},
$$
such that
\begin{enumerate}
\item[{\rm (1)}] each $E_{ij}$ is  an exponential equation over $H_{\lambda}$ for some $\lambda\in \{1,\dots,m\}$ or a trivial equation of kind $g_{ij}=1$, where $g_{ij}$ is an element of $G$,

\item[{\rm (2)}] for any $i=1,\dots, k$, the sets of variables of
$E_{i,j_1}$ and $E_{i,j_2}$ are disjoint if $j_1\neq j_2$,

\item[{\rm (3)}] $E$ is solvable if and only if $\Phi$ is solvable.\\
Moreover, any solution of $\Phi$ can be algorithmically extended to a solution of~$E$.

\end{enumerate}
}


{\it Proof.} Consider the exponential equation $E$, which is
$$
a_1g_1^{x_1}a_2g_2^{x_2}\dots a_ng_n^{x_n}=1\eqno{(8.3)}
$$
with $a_1,\dots,a_n,g_1,\dots,g_n\in G$.
Let $A_{par}, A_{fin}, A_{lox},$ be the subsets of $\{g_1,\dots,g_n\}$ consisting of parabolic elements,
non-parabolic elements of finite order, and loxodromic elements, respectively. We have
$$
\{g_1,\dots,g_n\}=A_{par}\sqcup A_{fin}\sqcup A_{lox}.
$$


If the equation $E$ has a solution then, by Theorem~A, there exists a solution $(k_1,\dots,k_n)$
with
$$
|k_j|\leqslant \Bigl(n^2+\overset{n}{\underset{i=1}\sum}\,|a_i|_{X\cup \mathcal{H}}+\overset{n}{\underset{i=1}{\sum}}\, |g_i|_{X\cup \mathcal{H}}\Bigr)\cdot M
$$
for all $g_j\in A_{lox}$. Hence, the solvability of $E$ is equivalent to the
solvability of a finite disjunction of equations of type (8.3) with $A_{lox}=\emptyset$.
Therefore, we assume that $A_{lox}=\emptyset$. For elements $g_j\in A_{fin}$, it is sufficient to look for solutions
with $k_j\in \{0,1,\dots,m_j-1\}$, where $m_j$ is the order of $g_j$.
Therefore we may additionally assume that $A_{fin}=\emptyset$.
Thus, we have reduced to the case where all elements $g_i$ are parabolic.
For any parabolic $g_i$, there exists $h_i\in G$ such that $h_i^{-1}g_ih_i\in H_{\lambda(i)}$ for some $\lambda(i)\in \{1,\dots,m\}$. This reduces the problem to Theorem B,
which gives the desired $\Phi$.
\hfill $\Box$



\def\refname{REFERENCES}

\end{document}